\documentclass[reqno]{amsart}
\usepackage{amssymb}
\usepackage{graphicx} 
\usepackage{amsmath}
\usepackage{amsthm}
\usepackage{amsfonts}
\usepackage{amstext}
\providecommand{\U}[1]{\protect\rule{.1in}{.1in}}
\theoremstyle{plain}
\newtheorem{theorem}{Theorem}[section]

\newtheorem{cor}[theorem]{Corollary}

\newtheorem{lem}[theorem]{Lemma}

\newtheorem{prop}[theorem]{Proposition}

\theoremstyle{definition}

\newtheorem{rem}[theorem]{Remark}
\numberwithin{equation}{section}
\numberwithin{theorem}{section}
\ifx\pdfoutput\relax\let\pdfoutput=\undefined\fi
\newcount\msipdfoutput
\ifx\pdfoutput\undefined\else
\ifcase\pdfoutput\else
\ifx\paperwidth\undefined\else
\ifdim\paperheight=0pt\relax\else\pdfpageheight\paperheight\fi
\ifdim\paperwidth=0pt\relax\else\pdfpagewidth\paperwidth\fi
\fi\fi\fi
\begin{document}
\title[Loop type subcontinua for indefinite concave-convex problems]{Loop type subcontinua of positive solutions for indefinite concave-convex problems}

\author{Uriel Kaufmann}
\address{U. Kaufmann \newline FaMAF, Universidad Nacional de C\'{o}rdoba, (5000) C\'{o}rdoba, Argentina}
\email{\tt kaufmann@mate.uncor.edu}

\author{Humberto Ramos Quoirin}
\address{H. Ramos Quoirin \newline Universidad de Santiago de Chile, Casilla 307, Correo 2, Santiago, Chile}
\email{\tt humberto.ramos@usach.cl}

\author{Kenichiro Umezu}
\address{K. Umezu \newline Department of Mathematics, Faculty of Education, 
Ibaraki University, Mito 310-8512, Japan} 
\email{\tt kenichiro.umezu.math@vc.ibaraki.ac.jp}


%{Concave-convex elliptic problem; Indefinite weight; Positive solution; Loop type subcontinuum; Bifurcation; Sub and supersolutions}

\thanks{U. Kaufmann was partially supported by Secyt-UNC 30720150100019CB} 
\thanks{H. Ramos Quoirin was supported by FONDECYT grants 1161635, 1171532 and 1171691}
\thanks{K. Umezu was supported by JSPS KAKENHI Grant Numbers 15K04945 and 18K03353}

%\date{\today}

\maketitle
\begin{abstract} 
We establish the existence of {\it loop type subcontinua} of nonnegative solutions for a class of concave-convex type elliptic equations with indefinite weights, under Dirichlet and Neumann boundary conditions. Our approach depends on local and global bifurcation analysis from the zero solution in a non-regular setting, since the nonlinearities considered are not differentiable at zero, so that the standard bifurcation theory does not apply.  To overcome this difficulty, we combine a regularization scheme with  {\it a priori} bounds, and Whyburn's topological method. Furthermore, via a continuity argument we prove a positivity property for subcontinua of nonnegative solutions. These results are based on a positivity theorem for the associated concave problem proved in \cite{KRQU}, and extend previous results established in the powerlike case. 
\end{abstract}

%------         start on 26Feb2017 
%------         start on 5Jun2018 for revision

\section{Introduction and main results} 
Let $\Omega\subset\mathbb{R}^{N}$ ($N\geq1$) be a bounded domain with smooth
boundary $\partial\Omega$. In this paper, we consider nonnegative solutions of the problem  
%\mpfns{modified, agree} 
%\mpfnsb{For tech.comment 8}
\[%
\begin{cases}
-\Delta u=\lambda a(x)f(u)+ b(x)g(u) & \mbox{ in }\Omega,\\
\mathcal{B}u=0 & \mbox{ on }\partial\Omega,
\end{cases}
\leqno{(P_{\mathcal{B}})}
\]
where:
\begin{itemize}
\item $\Delta$ is the usual Laplacian in $\mathbb{R}^{N}$; 
\item $\mathcal{B}u:= u$ (Dirichlet) or $\mathcal{B}u:=
\frac{\partial u}{\partial \mathbf{n}}$ (Neumann),
where $\mathbf{n}$ is the outward unit normal to $\partial \Omega$; 
\item $\lambda\in\mathbb{R}$ is a bifurcation parameter; 
\item $a, b \in C(\overline{\Omega})$ are such that 
%\mpfns{modified, agree}
$a$ changes sign in  $\Omega$ and $b(x_0)>0$ for some  $x_0 \in \Omega$;

\item $f,g:[0,\infty)\rightarrow \mathbb{R}$ are continuous functions with $f(0)=g(0)=0$. 
\end{itemize}
It follows that $(P_{\mathcal{B}})$ possesses the trivial line 
%\mpfnsb{For tech.comment 1}
$(\lambda,0)$ of zero solutions. 
The prototype of $f, g$ to be considered in this paper 
is %\mpfnsb{For tech.comment 4}
\begin{align} \label{pt}
f(s)=s^{q}, \quad g(s)=s^{p}, \quad\text{with }0<q<1<p,
\end{align}  so that the nonlinearity $\lambda f(s) + g(s)$ has a concave-convex nature. More precisely, 
we assume that $f \in C^{1}((0,\infty))$ with $f(s)>0$ for  $s>0$ and $g \in C^{1}([0, \infty))$ satisfy
\begin{align}
&  \lim_{s\rightarrow0^{+}}\frac{f(s)}{s}=\infty,\label{f/sinfty}\\
&  \lim_{s\rightarrow0^{+}}\frac{g(s)}{s}=0,\label{H_1g}\\
&  \lim_{s\rightarrow\infty}\frac{f(s)}{s}=0,\label{H_2f}\\
&  \lim_{s\rightarrow\infty}\frac{g(s)}{s}=\infty. \label{H_2g}%
\end{align}

Let $r>N$. A function $u\in W^{2,r}(\Omega)$ (and consequently, $u\in
C^{1}(\overline{\Omega})$) is said to be a \textit{nonnegative
solution} of $(P_{\mathcal{B}})$ if $u\geq0$ in $\Omega$, $u$ satisfies 
%\mpfns{paragraph modified, agree}
the equation pointwisely %\mpfns{modified}
a.e. in $\Omega$, and $\mathcal{B}u=0$ on $\partial\Omega$. If, in addition, $u$ satisfies \[ \begin{cases} u>0 \text{ in } \Omega \text{ and } \frac{\partial u}{\partial\mathbf{n}}<0 \text{ on } \partial \Omega & \text{ if } \mathcal{B}u=u,\\
 u>0 \text{ on } \overline{\Omega} & \text{ if } \mathcal{B}u=\frac{\partial u}
{\partial \mathbf{n}}, \end{cases} \] then we write $u\gg 0$. In this case $u$ lies in the interior of the positive cone $\{u\in V:u\geq0\}$, where
\begin{align*}
V:= \left\{ \begin{array}{ll}
C^1_0(\overline{\Omega}):=\left\{  u\in C^{1}(\overline{\Omega}):u=0\text{ on }\partial\Omega\right\} & \mbox { if }\ \ \mathcal{B}u=u, \\
C^1(\overline{\Omega}) & \mbox{ if }\ \ \mathcal{B}u=\frac{\partial u}
{\partial \mathbf{n}}.
\end{array}\right.
\end{align*}
A %\mpfns{added}
nonnegative solution $u$ of $(P_{\mathcal{B}})$ is called \textit{positive} if $u\gg 0$. 

%Let us remark that, due to the fact that $a$ changes sign and since $0<q<1$, there can exist solutions $u$ of $(P_{\mathcal{B}})$ (under either Dirichlet or Neumann boundary condition) such that $u>0$ in $\Omega$ but not satisfying that $u\gg 0$ (for concrete examples we refer to \cite{KRQU,KRQU2,KRQU3}).

Our first goal is to establish, under certain conditions on $a$ and $b$,  
the existence of  loop type \textit{subcontinua} $\{ (\lambda, u) \}$ 
(i.e., nonempty, closed and connected subsets in $\mathbb{R} \times V$)
composed by $(0,0)$ and nontrivial nonnegative solutions
$(\lambda, u)$ of $(P_{\mathcal{B}})$.
We shall prove the existence of a \textit{loop type} subcontinuum $\mathcal{C}_0$ such that %\mpfnsb{For tech.comment 2}
\begin{align} \label{single}
\mathcal{C}_{0} \cap \{ (\lambda,0) : \lambda \in \mathbb{R} \} = \{ (0,0)\}.
\end{align} 

It should be emphasized that, in general, one can not deduce that nontrivial nonnegative solutions of $(P_{\mathcal{B}})$ satisfy $u \gg 0$, since 
%\mpfns{modified, agree}
the strong maximum principle does not apply. This is due to the fact that
$a(x)f(\cdot)$ does not satisfy the slope condition \cite[p.623]{Am76}, 
see 
%\eqref{scd} below. 
%\marginpar{At some places we write Theorem ??(i) and at others Theorem ?? (i); and the same with Remarks, Figs, etc; I think we should unify the notation, both are ok for me.}
%\mpfns{you mean spacing, right ?}
Remark \ref{rem:H3}(ii) below. As a matter of fact, $(P_{\mathcal{B}})$ may have solutions $u$ satisfying
%(under either Dirichlet or Neumann boundary condition) 
 $u>0$ in $\Omega$ but not $u\gg 0$ (for concrete examples one may argue as in the proof of \cite[Proposition 2.9]{KRQU} after a slight modification).
In view of this difficulty, our second purpose  is to show 
that nontrivial solutions lying on 
$\mathcal{C}_0$ satisfy $u \gg 0$.

Bounded subcontinua of positive solutions for indefinite superlinear equations of the form
\begin{align*} 
-\Delta u = \lambda a(x)u + b(x)u^p \quad \mbox{ in } \ \Omega, 
\end{align*}
with $\Omega$ bounded (under different boundary conditions) or $\Omega = \mathbb{R}^N$, have been studied by several authors,  see e.g. \cite{B07,Ca04,CD96,CLGMM04,LGS04,CLGMM05,LGMM05,LG13}. According to %\mpfnsb{For tech.comment 3}
\cite{CLGMM04,CLGMM05,LGMM05}, a bounded subcontinuum linking two different points 
on $(\lambda,0)$ is called a \textit{mushroom}, one that meets a single point on $(\lambda,0)$ is called a \textit{loop}, and one that does not touch $(\lambda,0)$ is called an \textit{isola}. 
Cingolani and G\'amez studied both the Dirichlet condition case and the case $\Omega = \mathbb{R}^N$, proving the existence of mushrooms \cite[Theorems 4.4 and 5.5]{CD96}. 
Cano-Casanova considered a mixed boundary condition (with a second order uniformly strongly elliptic operator), and proved the existence of a mushroom \cite[Theorem 1.4]{Ca04}. 
L\'opez-G\'omez and Molina-Meyer dealt with the Dirichlet condition, and established existence results for a mushroom, a loop and an isola in three cases, respectively \cite[Theorems 3.1, 5.1 and 5.2]{LGMM05}. 
In the case of Neumann boundary conditions, Brown proved the existence of a mushroom and a loop in two situations, respectively \cite[Sections 2 and 5]{B07}. 
Finally, 
we refer to \cite[Section 3]{Um12} for the existence of a mushroom 
of positive solutions for a  semilinear elliptic problem with a logistic nonlinearity and an indefinite weight, coupled with a nonlinear boundary condition. 
Let us emphasize that all the previous works hold in the regular case, i.e.,  when the nonlinearity considered is $C^1$ at $u=0$, so that  
the general theory on local and global bifurcation from simple eigenvalues can be directly applied. 

Regarding existence results for positive solutions of concave-convex problems, a large number of works have been devoted to $(P_{\mathcal{B}})$ in the `definite case' (i.e. with $a \geq 0$, $a\not \equiv 0$) since the classical work of Ambrosetti, Brezis and Cerami \cite{ABC94}, which treats the model case \eqref{pt} 
with $a=b \equiv 1$ and %\mpfns{modified}
$p\leq \frac{N+2}{N-2}$ under the Dirichlet boundary condition. In \cite{ABC94} it is proved that $(P_{\mathcal{B}})$ has two positive
solutions for $\lambda>0$ sufficiently small. This result was extended by De Figueiredo, Gossez and Ubilla \cite{DGU2} to the non-powerlike case, with $a\geq0$. In addition, in 
\cite{DGU1}, the authors allowed $a$ to change sign and proved the existence
of two nontrivial nonnegative solutions of $(P_{\mathcal{B}})$ for $\lambda>0$ small. We refer to \cite{RQUIsrael} for a discussion on concave-convex problems under the Neumann boundary condition.

To the best of our knowledge, besides \cite{KRQU, RQUIsrael} there are no 
works providing the existence of solutions that are positive in $\Omega$ for 
%\mpfns{modified}
\textit{indefinite} concave-convex problems (i.e., with $a$ changing sign). In \cite{KRQU, KRQU2, KRQU3} we first established a positivity property for $(P_{\mathcal{B}})$ in the powerlike and concave case, i.e. with $f(s)=s^q$ and $b\equiv 0$. Thanks to these results, we obtained a positivity result for $(P_{\mathcal{B}})$ with $f(s)=s^q$ and $b \equiv 1$ (see \cite[Section 4]{KRQU}). Finally,
let us mention that in the model case \eqref{pt}, the existence of a loop type 
subcontinuum of nonnegative solutions for the Neumann case was obtained by means of a bifurcation approach in \cite{RQUIsrael}. Furthermore, the asymptotic profile of 
 nonnegative solutions as $\lambda \to 0^+$ enables one to deduce in some cases their positivity for $\lambda>0$ small, cf. \cite[Corollary 1.3]{RQUIsrael}.

For our first purpose, 
we assume that there exist two balls $B, B' \Subset \Omega$ 
and constants $a_0, a_0', b_0, b_0' > 0$ such that 
\begin{align} \label{ab:posi} 
\begin{cases}
a(x) \geq a_0 \ \mbox{ and } \ b(x) \geq b_0 & \mbox{ in } B, \\
-a(x) \geq a_0' \ \mbox{ and } \ b(x) \geq b_0' & \mbox{ in } B'.
\end{cases}
\end{align}
Let 
%\mpfns{modified, modified again} 
%\mpfnsb{For tech.comment 7}
$\psi \in C(\overline{\Omega})$ be such that \begin{align} \label{psi>0}
\Omega^\psi_+:=\{ x \in \Omega : \psi(x) > 0 \} \not =\emptyset.  
\end{align}
Then, we introduce the condition
\[
\mbox{$\Omega^\psi_+$ consists of a finite number of connected components of $\Omega$.}	\eqno{(H_\psi)}
\]
We shall assume this condition for $\psi=a$ and $\psi=-a$. 

%\mpfnsb{For tech.comment 4} %\mpfns{modified}
Motivated by the model case \eqref{pt}, we assume that
\begin{equation}
\lim_{s\rightarrow0^{+}}s^{1-q}f'(s)=:f_{0}\in(0,\infty
)\quad\mbox{for some $q\in (0,1)$.}\label{H_3}%
\end{equation}
We will see that under this condition $f$  
%\mpfns{modified}
behaves like $\frac{f_0}{q}s^q$ when $s\to 0^+$, and satisfies the slope condition, see Remark \ref{rem:H3}. 
In addition, the following \textit{strong} concavity (respect.\ 
convexity) condition on $f$ (respect.\ $g$) 
shall be used:
\begin{align}
& \left( \frac{f(s)}{s^q} \right)' \leq 0 \ \ \mbox{ for } s>0, \label{f:cave} \\
& \left( \frac{g(s)}{s} \right)' > 0 \ \ \mbox{ for } s>0,  \label{g:vex} 
\end{align}
where $q \in (0,1)$ is given by \eqref{H_3}. 
We introduce now the Gidas-Spruck condition \cite[Theorem 1.1]{GS81}, which is stronger than \eqref{H_2g}: 
\begin{align}
0<\lim_{s\to\infty} \frac{g(s)}{s^{p}} <\infty \quad 
\mbox{for some $p>1$, where $p<\frac{N+2}{N-2}$ if $N>2$}.  
\label{H_2'} 
\end{align}
%\mpfns{I think this sentence is not necessary, so we could remove it.}\mpfns{Modified. We use $(\mathcal{H}_b)$ in both thms.} 
Finally, %\mpfnsb{For tech.comment 7}
we shall use the condition $(\mathcal{H}_b)$, 
which will be precisely stated in Remark \ref{rem:H_b} 
and goes back to Amann and L\'opez-G\'omez  
\cite{ALG98}. 

For our second purpose, we focus on the case 
$f(s)=s^q$, $q\in (0,1)$.  
In association with the sublinear problem %\mpfns{modified}
\begin{align} \label{prb:aqD}
\begin{cases}
-\Delta u = a(x) u^{q} & \mbox{ in } \Omega,\\
\mathcal{B}u = 0 & \mbox{ on } \partial\Omega,
\end{cases}
\end{align}
we introduce the set
\begin{align} \label{defA}
\mathcal{A}_{\mathcal{B}}^{a} := \{ q \in(0,1) :u\gg 0 \ \ 
\mbox{for any nontrivial nonnegative solution $u$ of \eqref{prb:aqD}} \}.
\end{align}
We know \cite[Corollary 1.5 and Theorem 1.9]{KRQU} that 
under $(H_\psi)$ with $\psi = a$, there exists $q_{a} \in[0,1)$ 
such that $\mathcal{A}_{\mathcal{B}}^{a}=(q_{a}, 1)$, assuming additionally 
$\int_\Omega a < 0$ if $\mathcal{B}u=\frac{\partial u}{\partial \mathbf{n}}$.
Let us point out that the condition $\int_{\Omega}a<0$ is necessary and sufficient for the existence of a positive solution %$u\gg 0$
of \eqref{prb:aqD} with $\mathcal{B}u=\frac{\partial u}
{\partial \mathbf{n}}$, for some $q\in(0,1)$, 
see \cite[Corollary 1.3]{KRQU2}.

We are now in position to state our main results. 
First, we deal with the Dirichlet problem.

%     Main Theorem D

\begin{figure}[h] 
\centerline{
\includegraphics[scale=0.155]{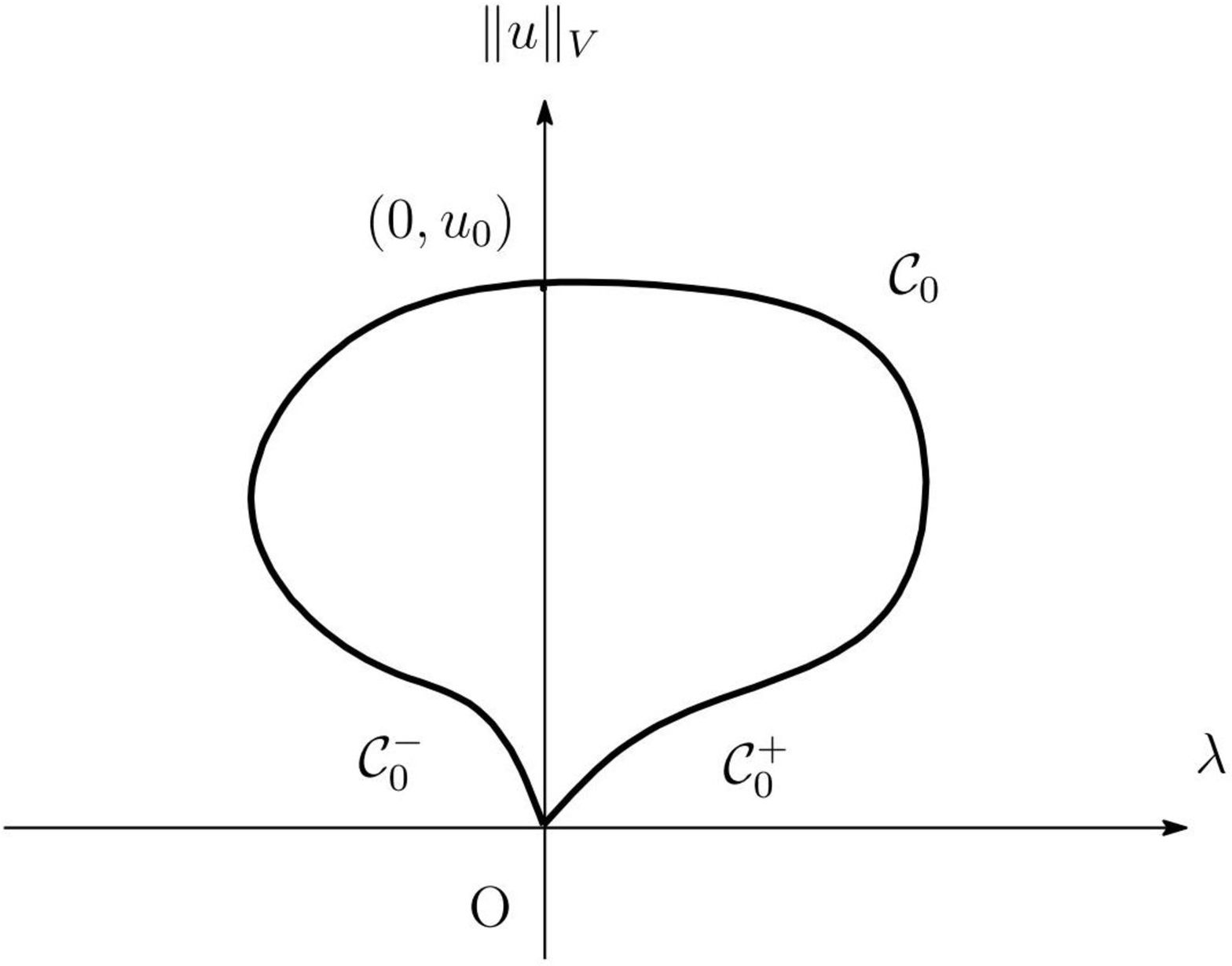} \hskip0.35cm
\includegraphics[scale=0.155] {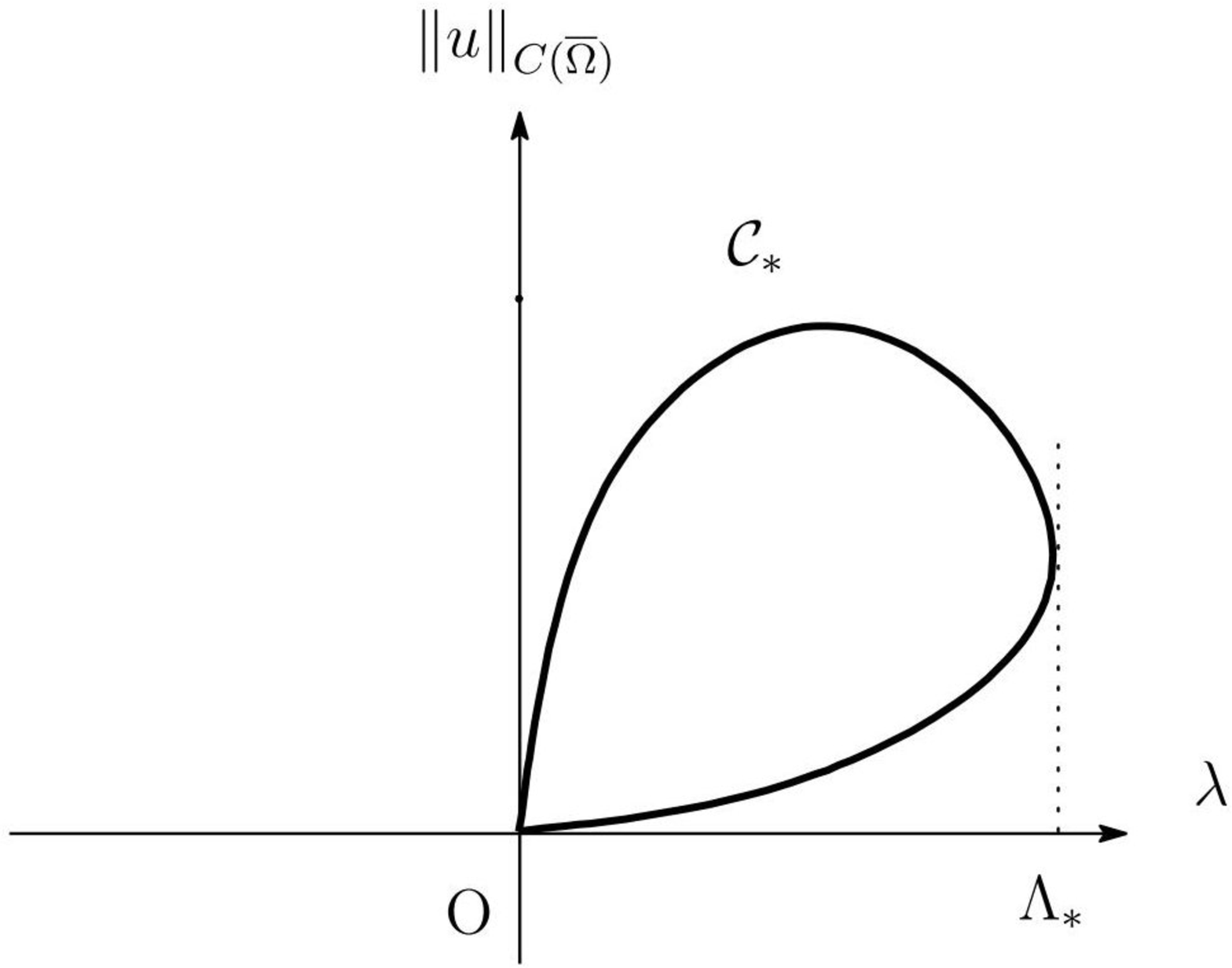}
} \centerline{(i) \hskip5.0cm (ii) }\caption{ The loop type subcontinua $\mathcal{C}_{0}$ and $\mathcal{C}_{\ast}$.} %
\label{fig:double}%
\end{figure}

\begin{theorem}
\label{resu01} Under 
%\mpfns{modified}
%\mpfnsb{For comment 3}
$\mathcal{B}u=u$, we assume  \eqref{H_1g}, \eqref{H_2f},  \eqref{ab:posi}, \eqref{H_3}, \eqref{H_2'},
%\mpfns{one $(\mathcal{H}_b)$ deleted. this condition is assumed in case (b).} %$(\mathcal{H}_b)$, 
and 
%\mpfnsb{modified}
$(H_\psi)$ with $\psi = \pm a$. In addition, suppose either 
\begin{enumerate}
\item[(a)] $b>0$ on $\overline{\Omega}$, and $p< \frac{N+2}{N-2}$ if $N>2$, or 
\item[(b)] \eqref{f:cave}, \eqref{g:vex}, and $(\mathcal{H}_b)$. 
%\mpfns{what is $(H_b)$?, clarified, I think}
\end{enumerate}
Then, the following two assertions hold:
\begin{enumerate}
\item 
$(P_{\mathcal{B}})$ admits a loop type \textrm{subcontinuum} 
$\mathcal{C}_{0}$ (i.e., a nonempty,
closed and connected subset in $\mathbb{R}\times C^1_0(\overline{\Omega})$) 
of nonnegative solutions which satisfies \eqref{single}. 
Moreover, we have the following properties, see 
Figure \ref{fig:double}(i):
\begin{enumerate}
\item $(0,u_{0}) \in \mathcal{C}_{0}$ for some positive solution $u_0$ of $(P_{\mathcal{B}})$ with $\lambda=0$. 
\item There exists $\delta> 0$ such that $\mathcal{C}_{0}$ does not 
contain any positive solution $u$ of $(P_{\mathcal{B}})$ with $\lambda = 0$ satisfying $\Vert u \Vert_{C(\overline{\Omega})} \leq\delta$. 
\item $\mathcal{C}_{0}$ contains closed connected sets 
$\mathcal{C}^{\pm}_{0}$ such that $\{ (0,0) \} \subsetneq 
\mathcal{C}^{\pm}_{0}$, and 
if $(\lambda, u) \in\mathcal{C}^{\pm}_{0} 
\setminus\{ (0,0)\}$, then $\lambda\gtrless0$, i.e., $\mathcal{C}%
_{0}$ bifurcates both subcritically and supercritically at $(0,0)$. 
\newline
\end{enumerate}

\item Let $f(s)=s^q$, $q\in (0,1)$.  
Assume that $b\geq 0$, and additionally that 
\begin{align}
\label{Hg:posi} g(s) \geq0 \ \ \mbox{ for } s> 0 
\end{align}
when condition (a) holds. If $q \in \mathcal{A}_{\mathcal{B}}^{a}\cap\mathcal{A}_{\mathcal{B}}^{-a}$, then $u \gg 0$ for any $(\lambda, u) \in\mathcal{C}_{0}
\setminus\{ (0,0) \}$.  In particular, the component of nontrivial nonnegative
solutions of $(P_{\mathcal{B}})$ including $\mathcal{C}_{0}\setminus 
\{ (0,0)\}$ is bounded. \newline
\end{enumerate}
%\end{enumerate}
\end{theorem}

%\mpfns{modified the caption, agree} 
%It \mpfns{maybe this sentence could be removed, since it's not so important. Agree. \textbf{It would be enough that this is mentioned in the proof.}}\mpfns{old Remark 1.7}should be remarked that from Remark \ref{rem:cc} (iii) below, we see that \eqref{Hg:posi} holds in case (b) of Theorem \ref{resu01}. 

Next we consider the Neumann problem under the condition 
%\mpfns{I think we should explain why we're assuming this condition. Agree. Removed $(x)$ in 1.16} 
\begin{align} \label{ib<0}
\int_\Omega b < 0. 
\end{align}
In this case, %\mpfns{modified} 
%and complement the results of \cite{RQUpre}, where the case $\int_\Omega b \geq 0$ was treated (under \eqref{pt}), 
we shall obtain a loop type subcontinuum of 
nonnegative solutions with the same nature the Dirichlet case admits, 
as in Figure \ref{fig:double}(i) (for the case $\int_\Omega b\geq 0$ we refer to 
Remark \ref{rem:ibgeq0} below). 
To this end, we need the following decay 
and positivity condition for $g$, which is stronger 
than \eqref{H_1g}: 
\begin{equation}
\lim_{s\rightarrow0^{+}}\frac{g(s)}{s^{\sigma}}=:
g_{0}\in(0,\infty) \quad 
\mbox{ for some $\sigma>1$, where $\sigma < \frac{2N}{N-2}$ 
if $N>2$}. 
\label{gssig} 
\end{equation}
Moreover, we are able to discuss the {\it positivity} of (nontrivial) 
nonnegative solutions for $(P_{\mathcal{B}})$ with \eqref{pt}, assuming 
\begin{align} \label{Nespe}
\left\{ 
\begin{array}{l}
p<\frac{N+1}{N-1} \quad \mbox{if $N>2$}, \\
a \in C^{\alpha}(\overline{\Omega}) \quad \mbox{for some 
$\alpha\in(0,1)$}, \\
\int_{\Omega}a <0, \\
b\equiv 1. 
\end{array}\right.
\end{align}
It is known \cite[Theorem 1]{RQUpre} that if $\Omega^a_+$ is connected, then $(P_{\mathcal{B}})$ possesses a 
loop type subcontinuum $\mathcal{C}_{\ast}$ in 
$\mathbb{R} \times C(\overline{\Omega})$ 
of nonnegative solutions which satisfies \eqref{single}. 
%Furthermore, we have the following properties, see Figure \ref{fig17_0105}: 
Furthermore, we have the following properties, see Figure \ref{fig:double}(ii): 
\begin{enumerate}

\item $\{ \lambda: (\lambda, u) \in\mathcal{C}_{\ast} \setminus\{ (0,0)\}
\} = (0, \Lambda_{\ast}]$ for some $\Lambda_{\ast} >0$. 

\item $\mathcal{C}_{\ast}$ possesses at least two nontrivial nonnegative
solutions for $\lambda> 0$ small enough. 
\end{enumerate}

Now, we state our main results for the Neumann problem, 
which are given in a similar way as in Theorem \ref{resu01}, 
and where condition \eqref{ib<0} provides us with a loop type 
subcontinuum bifurcating both subcritically and supercritically at $(0,0)$.

%             Main Theorem N

\begin{theorem} \label{resu02} 
Under %\mpfnsb{For comment 3}
$\mathcal{B}u=\frac{\partial u}{\partial \mathbf{n}}$, assume \eqref{ab:posi}, \eqref{H_3}, \eqref{f:cave}, \eqref{g:vex}, \eqref{H_2'}, \eqref{gssig}, $(H_\psi)$ with 
%\mpfnsb{modified}
$\psi = \pm a$, and 
%\mpfns{$\pm b$ deleted, $\mathcal{H}_b$ added}
$(\mathcal{H}_b)$. Then, 
%\marginpar{"," added}\mpfns{agree}
the following two assertions hold:
\begin{enumerate}
\item If \eqref{ib<0} holds, then $(P_{\mathcal{B}})$ admits a loop type 
subcontinuum 
$\mathcal{C}_{\ast}$ in $\mathbb{R}\times C^1(\overline{\Omega})$ of nonnegative 
solutions for which the same assertions in Theorem \ref{resu01}(i) 
hold true, see Figure \ref{fig:double}(i).

\item Assume \eqref{pt}, \eqref{Nespe} and the condition that $\Omega^a_+$ is connected. 
Let $\mathcal{C}_{\ast}$ be the loop type subcontinuum stated above. 
If $q \in \mathcal{A}_{\mathcal{B}}^{a}$, 
then the same conclusion in Theorem \ref{resu01}(ii) holds with
$\mathcal{C}_{0}$ replaced by $\mathcal{C}_{\ast}$.
\end{enumerate}
\end{theorem} 

%\mpfns{changed the place of remark and removed paragraph}
\begin{rem} \label{rem:ibgeq0} 
When $\mathcal{B}u=\frac{\partial u}{\partial \mathbf{n}}$ and 
$\int_\Omega b(x) \geq 0$, the existence of a  %bounded 
loop type subcontinuum 
of nonnegative solutions of $(P_{\mathcal{B}})$ has been established in the particular case $f(s)=s^q$ and $g(s)=s^p$ with $0<q<1<p$ (see \cite{RQUpre}, and as a particular case, see also \eqref{Nespe}). 
%\marginpar{I would remove "and also (1.18)", or clarify what we are meaning}\mpfns{modified}. 
In this case, 
%\mpfns{modified} 
although the loop type subcontinuum $\mathcal{C}_{0}$ satisfies \eqref{single}, $\mathcal{C}_{0}\setminus \{ (0,0)\}$ appears in $\lambda>0$. This means that $\mathcal{C}_{0}$ never meets the vertical line $\{ (0, u) : 0\not\equiv u \geq 0 \}$, see \cite[Lemma 6.8(1)]{RQUIsrael}. 
Thus, the approach used in the proof of Theorem \ref{resu02}(i) does not work for excluding the possibility that $\mathcal{C}_{0}=\{ (0,0)\}$, see the argument in Subsection \ref{subsec:prfi}. 
Let us mention that in \cite{RQUpre} the authors used a suitable rescaling technique (which strongly relies on the homogeneity of $f(s) = s^q$ and $g(s)=s^p$)  to exclude this possibility. \newline
\end{rem}

%It is worth pointing out that, \mpfnsb{For tech.comment 9}as a particular case of $f$ satisfying \eqref{H_3} and \eqref{f:cave} in Theorems \ref{resu01} and 
%\ref{resu02}, we give $f(s)=s^q/(1+s^{r})$ ($r\geq0$) and $s^qe^{-s}$, 
%based on Lemma \ref{lem:equiv}. We also remark that the case $\mathcal{B}u=\frac{\partial u}{\partial \mathbf{n}}$ and $\int_\Omega b\geq 0$ 
%will be mentioned in Remark \ref{rem:ibgeq0}.

\begin{rem}
Theorems %\mpfnsb{For tech.comment 7}
\ref{resu01} and \ref{resu02} can be extended to the case $a, b \in L^\infty (\Omega)$ except assertion (ii) in Theorem \ref{resu02}. This can be done if we formulate $(H_\psi)$ for $\psi \in L^\infty (\Omega)$ such that $\psi \not \equiv 0$, letting $\Omega^{\psi}_+$ be the largest open subset of $\Omega$ in which $\psi >0$ a.e., 
and assuming additionally %\mpfns{modified} 
$$ %(H_\psi) \ \
\begin{cases}
%& \Omega^\psi_+ =\displaystyle \bigcup_{i=1}^l \Omega_i, \mbox{ where $\Omega_i$ are open and connected for  $i=1,...,l$.  } 
& \mbox{$\psi$ is bounded away from zero on compact subsets of $\Omega_i$,} \\
& \mbox{for } i=1,...,l, \mbox{ and } |(\textrm{supp } \psi^+) \setminus \Omega^\psi_+| = 0. 
\end{cases}
$$
\newline
\end{rem}

\section{Preliminaries and Examples} \label{sec:preexa}
%\mpfns{changed name of section, no more section "examples"}
%\mpfns{modified}
We start this %\mpfnsb{For comment 2}
section with some remarks concerning some of our assumptions. 
\begin{rem}\label{rem:H3} Condition \eqref{H_3} implies:
\begin{enumerate}
\item by the L'Hospital rule, 
\begin{align}
\lim_{s\to 0^{+}} \frac{f(s)}{s^{q}} = \frac{f_0}{q}>0. \label{H_3'}%
\end{align}
In particular, since $f \in C^1((0,\infty))$ and $f(0)=0$, we can show that $f \in C^{\alpha}([0,s_{0}])$ for $\alpha \in (0,q]$ and $s_0>0$. 

\item $f$ satisfies the \textit{slope condition}, that is, for any $s_{0} > 0$,
there exists $M_{0} > 0$ such that
\begin{align}
\label{f:slope}\frac{f(s) - f(t)}{s-t} > - M_{0} \quad\mbox{ for } 0\leq t < s
\leq s_{0}.
\end{align}

\end{enumerate}

However, even under \eqref{H_3}, $a(x)f(\cdot)$ 
does \textit{not} satisfy the slope
condition for $x\in\Omega$ where $a(x)<0$, since $\displaystyle \lim_{s \to 0^+} f'(s)=\infty$. 
\end{rem}

%\mpfnsb{For tech.comment 7: old Rem1.2 moved here}
\begin{rem} \label{rem:cc}
\strut 
\begin{enumerate}
\item Since $f(s)>0$ for $s>0$, we note that if \eqref{f:cave} holds, then $f$
is concave for $s>0$, i.e.
\[
\left(  \frac{f(s)}{s}\right)  ^{\prime}<0\ \mbox{ for }s>0.
\]

\item It is easy to check that \eqref{f:cave} is stronger than \eqref{H_2f}.
This is a consequence of the fact that \eqref{f:cave} yields%
\[
0\leq \lim_{s\to \infty}\frac{f(s)}{s^q}<\infty. 
\]

\item Let us also note that \eqref{H_1g} and \eqref{g:vex} imply that $g(s)>0$
for $s>0$.
Indeed, assume first $g(s_{0})<0$ for some $s_{0}>0$, and set $\varepsilon
_{0}:=-g(s_{0})/s_{0}>0$. From \eqref{H_1g}, we infer that for some $s_{1}%
\in(0,s_{0})$,
\[
\frac{g(s_{0})}{s_{0}}=-\varepsilon_{0}<\frac{g(s_{1})}{s_{1}},
\]
which contradicts \eqref{g:vex}. Hence $g(s)\geq0$ for all
$s>0$. Next, assume that $g(s_{0})=0$ for some $s_{0}>0$. By \eqref{g:vex}, it
follows that $g(s)\not \equiv 0$ for $s\in(0,s_{0})$. This implies that
$g(s_{1})>0$ for some $s_{1}\in(0,s_{0})$. It follows that
\[
\frac{g(s_{1})}{s_{1}}>0=\frac{g(s_{0})}{s_{0}},
\]
which contradicts \eqref{g:vex} again, as desired.
\end{enumerate}
\end{rem}

\begin{rem} \label{rem:H_b}
We describe here
%\marginpar{"here" added}\mpfns{agree}
the explicit formula for the growth 
%\mpfns{put calligraphic in order to distinguish with $H_\psi$}
condition $\mathcal{H}_b$ of $b^+$ in a neighborhood of $\partial \Omega^b_+$ used in Theorems \ref{resu01} and \ref{resu02}, which originates from Amann and L\'opez-G\'omez \cite[Theorem 4.3]{ALG98}:
\[
(\mathcal{H}_{b})\ \
\begin{cases}
& \Omega^b_+ \ \mbox{ is a subdomain of }\Omega
\mbox{ with smooth boundary }\partial\Omega_{+}^{b}\mbox{ and either }\\
& \bullet\ \ \overline{\Omega_{+}^{b}}\subset\Omega
,\mbox{ and }b<0\mbox{ in }D_{b}:=\Omega\setminus\overline{\Omega_{+}%
^{b}},\mbox{ or }\\
& \bullet\ \ \Omega_{+}^{b}\supset\{x\in\Omega:d(x,\partial\Omega
)<\sigma\}\ \mbox{ for some }\sigma>0,\mbox{ and }b<0\mbox{ in }D_{b}.\\
& \mbox{In addition, }D_{b}\ \mbox{ is a subdomain of }\Omega
\mbox{ with smooth boundary, and}\\
&
\mbox{there exist $\gamma > 0$ and a function $\beta$ defined in a tubular neighborhood}\\
& U:=\{x\in\Omega_{+}^{b}:d(x,\partial\Omega_{+}^{b})<\sigma
\}\mbox{ of $\partial \Omega^b_+$ in $\Omega^b_+$, which is continuous,}\\
& \mbox{positive
and bounded away from zero, and satisfies}\\
& b^{+}(\cdot)=\beta(\cdot)d(\cdot,\partial\Omega_{+}^{b})^{\gamma}\text{ in
}U\text{ and }1<p<\min\left(  \frac{N+2}{N-2},\frac{N+1+\gamma}{N-1}\right)
\text{ if }N>2.
\end{cases}
\]%
\end{rem}

	%	% \begin{figure}[H] 
		% \begin{figure}[!htb]
		%\begin{center}
		%\includegraphics[scale=0.23]{02.eps} 
		%	\caption{Possible subdomains $\Omega^b_+$ satisfying $(H_b)$.} 
		%	\label{fig17_0216}
		%\end{center}		  
		 %\end{figure}

%\begin{rem} \label{rem:ibgeq0} 
%When $\mathcal{B}u=\frac{\partial u}{\partial \mathbf{n}}$ and 
%$\int_\Omega b(x) \geq 0$, the existence of a  %bounded 
%loop type subcontinuum 
%of nonnegative solutions of $(P_{\mathcal{B}})$ has been established in the particular case $f(s)=s^q$ and %%$g(s)=s^p$ with $0<q<1<p$ (see \cite{RQUpre}, and also \eqref{Nespe}). 
%In this case, \mpfns{modified} although the loop type subcontinuum $\mathcal{C}_{0}$ satisfies \eqref{single}, $\mathcal{C}_{0}\setminus \{ (0,0)\}$ appears in $\lambda>0$. This means that $\mathcal{C}_{0}$ never meets the vertical line $\{ (0, u) : 0\not\equiv u \geq 0 \}$, see \cite[Lemma 6.8(1)]{RQUIsrael}. 
%Thus, the approach used in the proof of Theorem \ref{resu02}(i) does not work for excluding the possibility that $\mathcal{C}_{0}=\{ (0,0)\}$, see the argument in Subsection \ref{subsec:prfi}. 
%Let us mention that in \cite{RQUpre} the authors used a suitable rescaling technique (which strongly relies on the homogeneity of $f(s) = s^q$ and $g(s)=s^p$)  to exclude this possibility. \newline
%\end{rem}

We conclude this section showing some examples of functions satisfying the previous conditions. We start with the following lemma, which characterizes the functions satisfying \eqref{H_3} and \eqref{f:cave}.

%\mpfnsb{modified}

%\mpfnsb{For tech.comment 9}

\begin{lem} \label{lem:equiv}
Let $f\in C(\left[  0,\infty\right)  )\cap C^{1}((0,\infty))$ with $f(0)=0$
and $f(s)>0$ for $s>0$. Then, the following two conditions are equivalent:
\begin{enumerate}
\item \eqref{H_3} and \eqref{f:cave} hold.
\item $f(s)=s^{q}h(s)$ for some $q\in(0,1)$ and $h\in C(\left[  0,\infty\right)  )\cap C^{1}((0,\infty))$ such that 
\[
h\text{ is nonincreasing,}\quad h(s)  >0 \text{ for } s \geq 0\quad\text{and\quad}%
\lim_{s\rightarrow0^{+}}sh^{\prime}\left(  s\right)  =0.
\]
\end{enumerate}
\end{lem}

\noindent\textit{Proof.} 
It is easy to see that $f$ as in condition (ii) fulfills \eqref{H_3} and \eqref{f:cave}. Conversely, if $f$ satisfies the aforementioned conditions, defining $h\left(  s\right)  :=s^{-q}f\left(  s\right)  $ for $s>0$ and $h\left(
0\right)  :=\lim\limits_{s\rightarrow0^{+}}h\left(  s\right)  $, it is also easy to check that $h$ has the desired properties. \qed

As particular cases, we mention $h(s)=\dfrac{1}{1+s^{r}}$
($r\geq0$) and $h(s)=e^{-s}$. \newline

%\mpfnsb{modified} 
We note that oscillatory cases are out of our scope. For instance, consider 
$h(s) = \sin \left( 
\frac{1}{s} \right)+2$. If we put $f(s):=s^q h(s)$ for $s>0$, and $f(0):=0$, 
then $f \in C([0,\infty)) \cap C^1((0,\infty))$ with $f(0)=0$ and $f(s)>0$ for $s>0$. Moreover, $f$ fulfills \eqref{f/sinfty} and \eqref{H_2f}, but \eqref{H_3'}, \eqref{f:slope} and \eqref{f:cave} do not hold. \newline

%\mpfnsb{modified} 
We now exhibit examples of $g$ satisfying \eqref{g:vex}, \eqref{H_2'}, and \eqref{gssig}. 

(a) We set%
\begin{align*}
g\left(  s\right)   &  :=s^{p}h\left(  s\right)  ,\quad\text{with\quad
}1<p<\frac{N+2}{N-2}\quad\text{if}\quad N>2,
\end{align*}
where $h \in  C^{1}\left(  \left[  0,\infty\right)  \right)$
 is
nondecreasing, bounded,
and satisfies one of the following conditions:
\begin{enumerate}
\item[(i)] $0=h\left(  0\right)  <h^{\prime}\left(  0\right)  .$

\item[(ii)] $0<h\left(  0\right)  .$
\end{enumerate}

Note that \eqref{gssig} holds if we choose $\sigma:=p+1$ in (i),
and  $\sigma:=p$ in (ii). An example for (i) is $h(s)=1-e^{-s}$, 
while $h(s)=\arctan\left(  s+ 1 \right)  $ is included in (ii). Another example  is given by $h(s)=\frac{s^r}{1+s^r}$, with $r=0$ or $r=1$. The case $r=1$ satisfies (i), whereas $r=0$ satisfies (ii). More generally, the function $g(s) = \frac{s^{p+r}}{1+s^r}\ (0\leq r<\frac{2N}{N-2}-p)$ fulfills \eqref{g:vex}, \eqref{H_2'}, and \eqref{gssig}. Indeed, we can take $\sigma := p+r$ for \eqref{gssig}. \newline

(b) Let $k\geq1$ and $1<p<\frac{N+2}{N-2}$ if $N>2$. We put 
\[
g(s):=s^{p}\frac{k+s}{1+s}.
\]
Then, $g$ satisfies \eqref{H_2'} and \eqref{gssig}. It is also
clear that \eqref{g:vex} holds when $k=1$. Meanwhile, when $k>1$, it satisfies
\eqref{g:vex} if additionally
\[
p>p_{1}(k),\ \mbox{ where }\ p_{1}(k):=\frac{2\sqrt{k}}{\sqrt{k}+1}.
\]
We note that $p_{1}(k)$ is increasing for $k>1$, and $p_{1}(k)\searrow1$ as
$k\rightarrow1^{+}$, whereas
$p_{1}(k)\nearrow2$ as $k\rightarrow\infty$. Let us finally observe that case
(b) is {\it not} included in any of the possibilities considered in case (a). 
Indeed, $h(s)=\frac{k+s}{1+s}$ is decreasing for $s\geq 0$. \newline

\section{Regularization schemes and transversality conditions}
Let %\mpfnsb{For tech.comment 4: the material in old page 4 moved here.}
us now  explain our approach to study bifurcation of nontrivial nonnegative solutions for $(P_{\mathcal{B}})$ from $(\lambda,0)$. From \eqref{f/sinfty},  we see that $f$ is not differentiable at $s=0$, so that we can not directly apply the usual bifurcation theory from simple eigenvalues to $(P_{\mathcal{B}})$. To overcome this difficulty, we proceed as in \cite{RQUIsrael,RQUpre}, `regularizing' $(P_{\mathcal{B}})$ at $u=0$, using $\varepsilon >0$. 
%\mpfnsb{For tech.comment 6}
We refer to \cite[Section 5]{LGMM05} for a similar approach introducing a new
parameter for a different regular problem.

We extend $g$ to $\mathbb{R}$ as a $C^{1}$ function 
and set $F:\mathbb{R}\rightarrow
\mathbb{R}$ by%
\begin{equation}
F(s):=\left\{
\begin{array}
[c]{ll}%
s^{1-q}f(s), & s\geq0,\\
\frac{f_{0}}{q}s, & s<0.
\end{array}
\right.  \label{def:F}%
\end{equation}
For $\varepsilon>0$ we shall study the auxiliary problem%
\[%
\begin{cases}
-\Delta u=\lambda a(x)(u+\varepsilon)^{q-1}F(u)+b(x)g(u) 
& \mbox{ in }\Omega,\\
\mathcal{B}u=0 & \mbox{ on }\partial\Omega.
\end{cases}
\leqno{    (P_{\mathcal{B},\varepsilon})   }
\]
Note that \eqref{H_3} and \eqref{H_3'} imply that $F\in C^{1}(\mathbb{R})$, $F(0)=0$ and
$F'(0)=\frac{f_0}{q}$, so that
\begin{equation}
s\mapsto(s+\varepsilon)^{q-1}F(s)\in 
C^{1}((-\varepsilon,\infty)).\label{s+epC1}
\end{equation}
Observe also that $(P_{\mathcal{B},0})$ corresponds 
to $(P_{\mathcal{B}})$, as far as nonnegative solutions are concerned.

Let us set
\[
h_{\lambda,\varepsilon}(x,s):=\lambda a(x)
(s+\varepsilon)^{q-1}F(s) + b(x)g(s).
\]
From 
\eqref{s+epC1} we see that $h(x,\cdot)$ 
satisfies the slope condition. 
Consequently, given a nontrivial nonnegative solution 
$u$ of $(P_{\mathcal{B},\varepsilon})$, we can choose $M>0$ 
such that $(-\Delta+M)u\geq0$ and
$\not \equiv 0$ in $\Omega$. Thus, by the strong 
maximum principle %\mpfns{modified}
and Hopf's lemma, we deduce $u\gg 0$, 
%, we deduce that $u>0$ in $\Omega$. Moreover, by, $\frac{\partial u}{\partial\mathbf{n}}<0$ on $u^{-1}(0) \cap \partial\Omega$, so that $u\gg 0$, 
%\mpfnsb{For tech.comment 5}
see \cite{GT01}, \cite[Theorem 7.10]{LG13}.

We shall then consider the linearized eigenvalue problem 
at $u=0$ for the regular problem $(P_{\mathcal{B},\varepsilon})$:  
\begin{equation}%
\begin{cases}
-\Delta\phi=\lambda a(x)\frac{f_{0}}{q}
\varepsilon^{q-1}\phi &
\mbox{ in }\Omega,\\
\mathcal{B}\phi=0 & \mbox{ on }\partial\Omega.
\end{cases}
\label{lep}%
\end{equation}
Since $a$ changes sign, \eqref{lep} has exactly two principal eigenvalues 
$\lambda_{1,\varepsilon}^{-}< 0<\lambda_{1,\varepsilon}^{+}$ 
(respect.\ $\lambda_{1,\varepsilon}^{-}= 0<\lambda_{1,\varepsilon}^{+}$) 
if 
\begin{align*}
\mathcal{B}u=u \quad \left(  \mbox{respect.} \ \ \mathcal{B}u=\frac{\partial u}{\partial \mathbf{n}}  \ \ \mbox{and} \ \ \int_\Omega a < 0   \right), 
\end{align*}
which are both simple, and furthermore, 
$(\lambda_{1,\varepsilon}^{\pm}, 0)$ satisfy 
the \textit{Crandall-Rabinowitz transversality condition}, 
%\mpfnsb{For tech.comment 5}
see \cite[Theorem 9.4]{LG13}. 

Thanks to the simplicity and transversality condition, 
the local bifurcation theory \cite[Theorem 1.7]{CR71} ensures the existence and uniqueness of positive solutions of $(P_{\mathcal{B},\varepsilon})$ bifurcating at $(\lambda_{1,\varepsilon}^{\pm},0)$. Moreover, the unilateral global bifurcation theory \cite[Theorem 6.4.3]{LG01} (see also \cite[Theorem 1.27]{Ra71}) ensures that $(P_{\mathcal{B},\varepsilon})$ possesses two  components $\mathcal{C}_{\varepsilon}^{\pm} = \{ (\lambda, u) \}$ (i.e., maximal, nonempty, closed and connected subsets in $\mathbb{R}\times V$) of nonnegative solutions 
emanating from  
$(\lambda_{1,\varepsilon}^{\pm},0)$, respectively (see 
Remark \ref{sec:appen02}). In addition,  
$\mathcal{C}_{\varepsilon}^{+}\setminus
\{(\lambda_{1,\varepsilon}^{\pm},0)\}$ and 
$\mathcal{C}_{\varepsilon}^{-}\setminus
\{(\lambda_{1,\varepsilon}^{\pm},0)\}$ 
consist of positive solutions. %(more precisely, solutions in $\mathcal{P}^\circ_{\mathcal{B}}$). 
This is due to  elliptic regularity and the fact (see \cite[Proposition 18.1]{Am76}) that $(P_{\mathcal{B},\varepsilon})$ has no bifurcating positive solutions from $(\lambda,0)$ at any $\lambda\neq\lambda_{1,\varepsilon}^{\pm}$.

Under some additional growth condition on $g$, we shall verify that $\mathcal{C}_{\varepsilon}^{\pm}$ are bounded in $\mathbb{R}\times V$ uniformly in $\varepsilon \in (0, 1]$, so that $\mathcal{C}_{\varepsilon}^{-}=\mathcal{C}_{\varepsilon}^{+}(:=\mathcal{C}_{\varepsilon})$ (i.e. $\mathcal{C}_{\varepsilon}$ is a mushroom). By simple computations, it can be shown easily that 
\begin{align} \label{peg0}
\lambda_{1,\varepsilon}^{\pm} \to 0 \quad\mbox{as} \ \ \varepsilon \to 0^+,
\end{align}
so that, passing to the limit as $\varepsilon\rightarrow0^{+}$, we shall observe by Whyburn's topological argument 
\cite[(9.12) Theorem]{Wh64} that 
\[
\mathcal{C}_{0}:=\limsup_{\varepsilon\rightarrow0^{+}}
\mathcal{C}_{\varepsilon}
\]
is a loop type subcontinuum which consists of nonnegative solutions of $(P_{\mathcal{B}})$ and satisfies \eqref{single}.

\begin{figure}[h]
\centerline{
\includegraphics[scale=0.17]{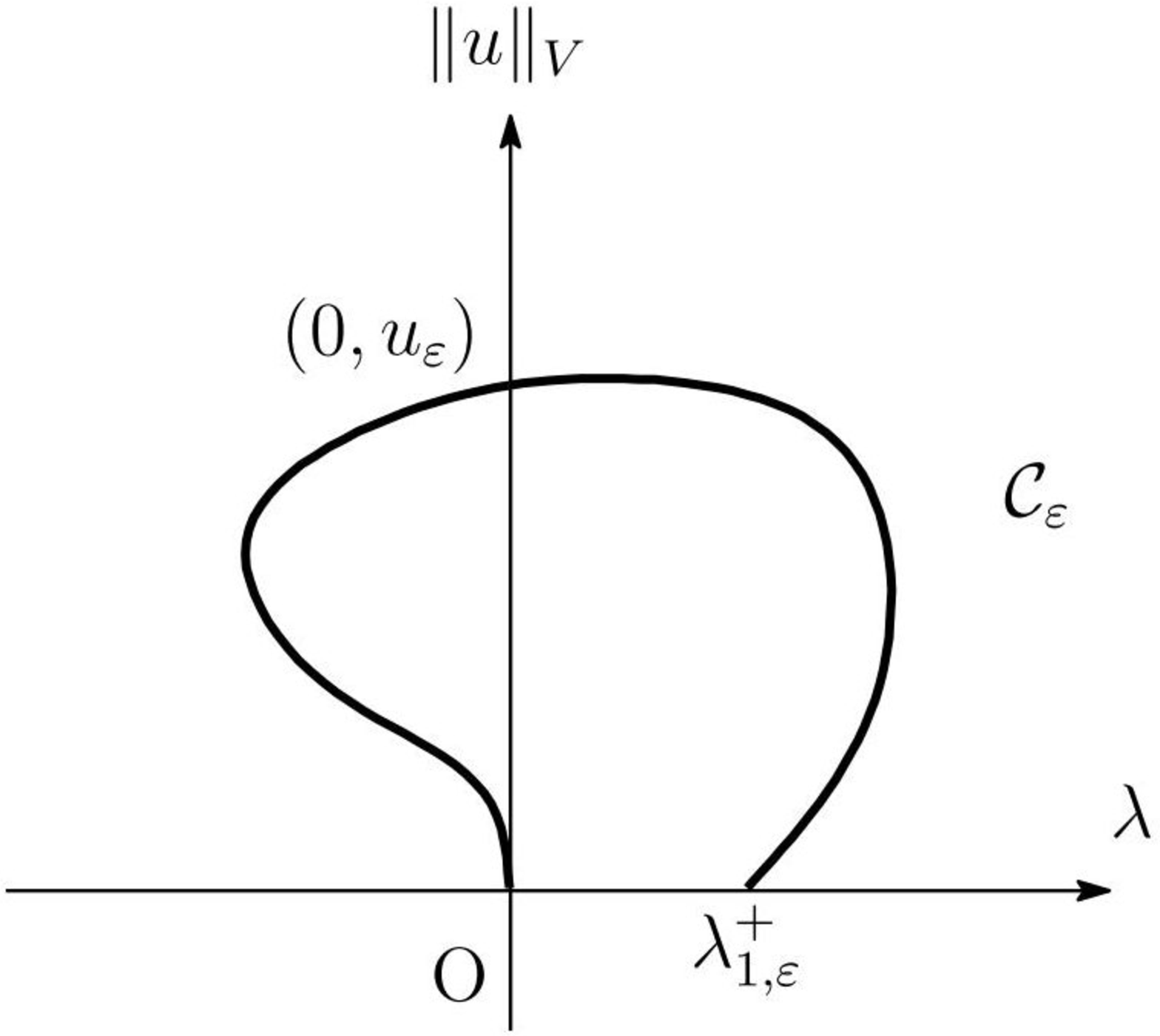} \hskip0.35cm
\includegraphics[scale=0.17] {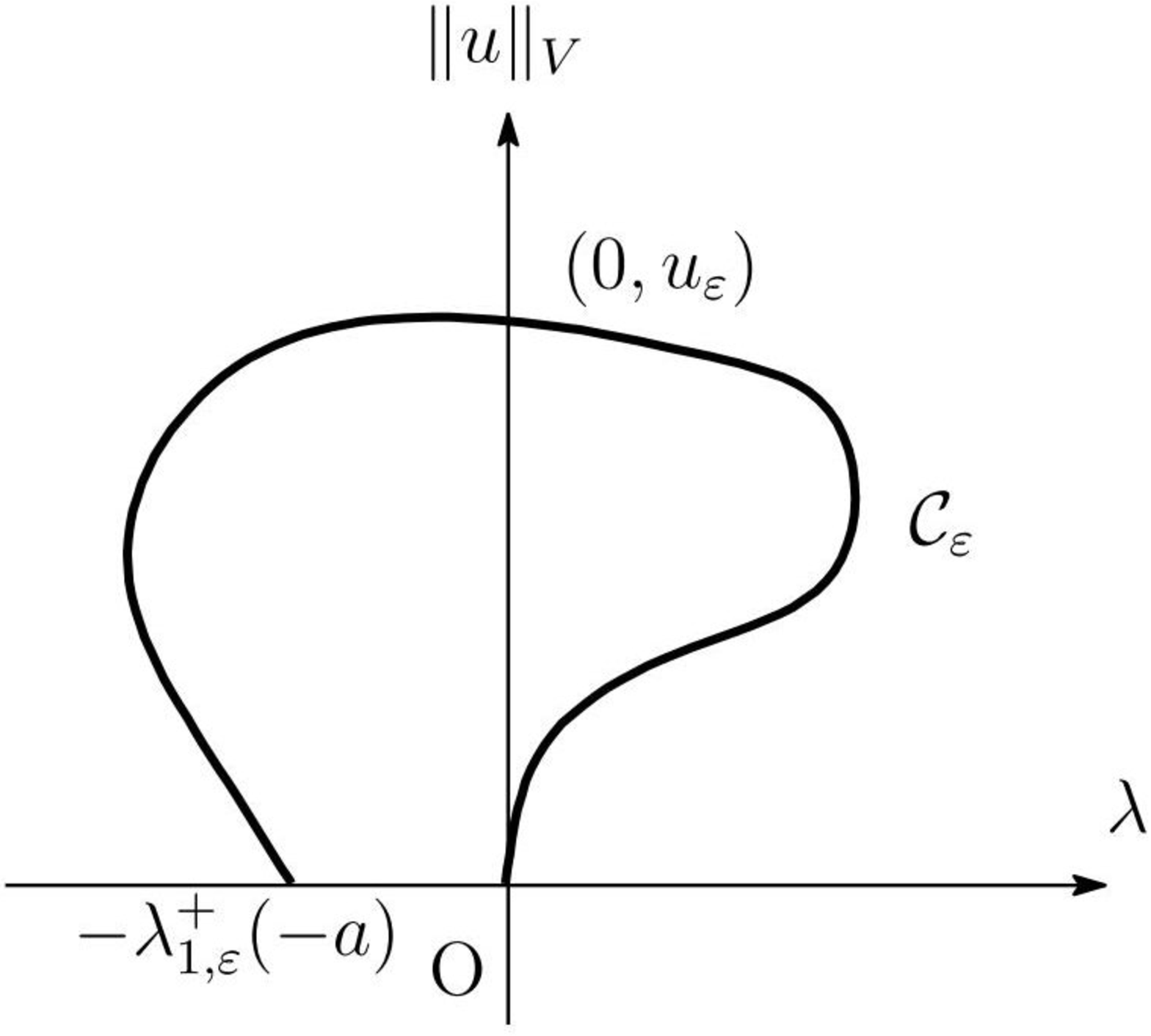}
} \centerline{(i) \hskip5.0cm (ii) }\caption{The bounded component $\mathcal{C}_{\varepsilon}$ for $(P_{\mathcal{B},\varepsilon})$ with $\mathcal{B}u=\frac{\partial u}{\partial \mathbf{n}}$. (i) Case $\int_\Omega a < 0$. (ii) Case $\int_\Omega a > 0$.}%
\label{fig:doublee}%
\end{figure}

\begin{rem} \label{case>0}
The critical case
\begin{align} \label{excepc}
\mathcal{B}u=\frac{\partial u}{\partial \mathbf{n}} \ \ \mbox{and} \ \ 
\int_\Omega a = 0
\end{align}
can be handled in a similar way. In this case, we replace $a$ by $a-\varepsilon$ for $\varepsilon > 0$ small in %\mpfns{modified}
$(P_{\mathcal{B},\varepsilon})$. Then, the above argument remains valid, since we can determine the asymptotic behavior \eqref{peg0} (see \cite[Lemma 6.6]{RQUIsrael} for the proof). In addition, we can reduce the case $\int_\Omega a > 0$ to the case $\int_\Omega a < 0$ under 
$\mathcal{B}u=\frac{\partial u}{\partial \mathbf{n}}$. 
Indeed, we only have to notice the symmetry property $\lambda a(x) = 
(-\lambda)(-a(x))$. The situation may be illustrated by Figures 
\ref{fig:doublee}(i) and (ii). 
\end{rem}

	%	% \begin{figure}[H] 
	%	 \begin{figure}[!htb]
	%	\begin{center}
	%	\includegraphics[scale=0.17]{mushrp.eps} 
	%		\caption{The bounded component $\mathcal{C}_{\varepsilon}$ for $(P_{\mathcal{B},\varepsilon})$ with $\mathcal{B}u=\frac{\partial u}{\partial \mathbf{n}}$ in the case $\int_\Omega a < 0$.} 
	%		\label{fig17_1010}
	%	\end{center}		  
	%	 \end{figure}

	%	% \begin{figure}[H] 
%		 \begin{figure}[!htb]
%		\begin{center}
%		\includegraphics[scale=0.16]{mushrn.eps} 
%			\caption{The bounded component $\mathcal{C}_{\varepsilon}$ for $(P_{\mathcal{B},\varepsilon})$ with $%\mathcal{B}u=\frac{\partial u}{\partial \mathbf{n}}$ in the case $\int_\Omega a > 0$.} 
	%		\label{fig17_1010b}
	%	\end{center}		  
	%	 \end{figure}

\begin{rem}\label{sec:appen02}
We %\mpfnsb{For tech.comment 14; old Appendix B}
shall show that the transversality condition allows us to apply the unilateral global bifurcation result \cite[Theorem 6.4.3]{LG01} to 
$(P_{\mathcal{B},\varepsilon})$ at $(\lambda_{1,\varepsilon}^{\pm}, 0)$. 
To this end, we reduce $(P_{\mathcal{B},\varepsilon})$ to an operator equation in $C(\overline{\Omega})$. Given $\xi \in L^r(\Omega)$, 
$r>N$, let $u \in W^{2,r}_{\mathcal{B}}(\Omega)$ be the unique solution of 
%\mpfns{modified: $(-\Delta + 1)$}
\begin{align} \label{pN}
\begin{cases}
(-\Delta +1) u = \xi (x) & \mbox{ in } \Omega, \\
\mathcal{B}u = 0 & \mbox{ on } \partial \Omega, 
\end{cases}
\end{align}
where 
\begin{align*}
W^{2,r}_{\mathcal{B}}(\Omega) := \left\{ \begin{array}{ll}
\{ u \in W^{2,r}(\Omega) : u=0 \ \mbox{on}\ \partial \Omega \} & \mbox{if} 
\ \mathcal{B}u = u, \\
W^{2,r}(\Omega) & \mbox{if} 
\ \mathcal{B}u = \frac{\partial u}{\partial \mathbf{n}}. 
\end{array} \right. 
\end{align*}
We introduce the solution operator $\mathcal{S} : L^r(\Omega) 
\to W^{2,r}_{\mathcal{B}}(\Omega)$ associated with \eqref{pN}, implying that $\mathcal{S}(\xi)=u$, which is bijective and homeomorphic. It follows that $\mathcal{S}: L^\infty (\Omega) \to W^{2,r}_{\mathcal{B}}(\Omega)$ is continuous, and moreover, $\mathcal{S}: L^\infty (\Omega) \to V$ is compact, since so is the embedding $W^{2,r}_{\mathcal{B}}(\Omega) \subset V$.  Thus, $(P_{\mathcal{B},\varepsilon})$ is reduced to 
\begin{align} \label{eq:cF}
\mathcal{F}(\lambda, u) := 
u - \mathcal{S}\left[ \left(  \lambda a \frac{f_0}{q}\varepsilon^{q-1} + 1 
\right) u + h_\lambda (x,u) \right] = 0 
\quad \mbox{ in } \quad C(\overline{\Omega}), 
\end{align}
where 
\[
h_\lambda(x,s) := \lambda a(x) \left\{ (s+\varepsilon)^{q-1}F(s) 
- \frac{f_0}{q}\varepsilon^{q-1}s \right\} +b(x)g(s). 
\]
Given 
$u \in C(\overline{\Omega})$, let
\[
\mathcal{A}_\lambda u:= \mathcal{S}\left[ \left(  \lambda a\frac{f_0}{q}\varepsilon^{q-1} + 1  \right)u \right], 
\]
and let 
\[
\mathrm{Ind}(0, \mathcal{A}_\lambda) := \deg (1-\mathcal{A}_\lambda, B_R) 
\]
be the fixed point index of $\mathcal{A}_\lambda$ at the 
origin for $\lambda \neq \lambda_{1,\varepsilon}^{\pm}$ but close to $\lambda_{1,\varepsilon}^{\pm}$, where $B_R$ is the ball with radius $R>0$ and centered at the origin. 
%Note that $\mathrm{Ind}(0, \mathcal{A}_\lambda)$ is well defined, since \cite[Proposition 18.1]{Am76} ensures that $u=0$ is the unique fixed point of $\mathcal{A}_\lambda$. 
%
%
According to \cite[Section 4.2, Theorem 5.6.2]{LG01},
%(see also \cite[p.419]{LGMM05}), 
the transversality condition at $(\lambda_{1,\varepsilon}^{\pm},0)$ implies that $\mathrm{Ind}(0, \mathcal{A}_\lambda)$ changes sign as $\lambda$ crosses $\lambda_{1,\varepsilon}^{\pm}$. Thus, \cite[Theorem 6.4.3]{LG01} applies, so that equation \eqref{eq:cF} (and so, $(P_{\mathcal{B},\varepsilon})$) possesses two components $\mathcal{C}_{\varepsilon}^{\pm} = \{ (\lambda, u) \}$ in $\mathbb{R}\times C(\overline{\Omega})$ of nonnegative solutions $u$ emanating from $(\lambda_{1,\varepsilon}^{\pm},0)$, respectively. 
%In addition, $\mathcal{C}_{\varepsilon}^{+}\setminus\{(\lambda_{1,\varepsilon}^{\pm},0)\}$ and $\mathcal{C}_{\varepsilon}^{-}\setminus\{(\lambda_{1,\varepsilon}^{\pm},0)\}$ consist of solutions in $\mathcal{P}^\circ_{\mathcal{B}}$, using the strong maximum principle and Hopf's lemma. 
%
Finally, we can verify that $\mathcal{C}_{\varepsilon}^{\pm}$ are also components in $\mathbb{R}\times V$, by elliptic regularity. 
%
%The closedness is clear. To obtain the connectivity, we assume by contradiction that there exist nonempty closed subsets $A, B$ of $\mathbb{R}\times V$ such that 
%	\[
%	\mathcal{C}_{\varepsilon}^{\pm} = A \cup B, \quad A \cap B = \emptyset. %
%	\]
%Then, we shall show that $A, B$ are nonempty and closed in $\mathbb{R}\times C(\overline{\Omega})$. Assume $(\lambda_n, u_n) \in A$ and $(\lambda_n, u_n) \to (\lambda_0, u_0)$ in $\mathbb{R}\times C(\overline{\Omega})$. From \eqref{eq:cF}, it follows that $\| u_n \|_{W^{2,r}_{\mathcal{B}}(\Omega)}$ is bounded. Using the compact embedding $W^{2,r}_{\mathcal{B}}(\Omega) \subset V$ again, we deduce that, up to a subsequence, $u_n \to u_0$ in $V$. By the assumption, we infer that $(\lambda_0, u_0) \in A$, as desired. The argument for the set $B$ is the same. However, this is contradictory to the connectivity in $\mathbb{R}\times C(\overline{\Omega})$. The connectivity in $\mathbb{R}\times V$ has been now obtained.%
\end{rem}

In the forthcoming sections, we will characterize the limiting behavior of $\mathcal{C}_{\varepsilon}$ as $\varepsilon\rightarrow0^{+}$ under the
%\marginpar{removed "certain"}\mpfns{agree}
conditions on $a,b,f$ and $g$ stated in Theorems \ref{resu01} and \ref{resu02}.

\section{A priori bounds} \label{sec:bounds}
In %\mpfnsb{For comment 4}
this section, we establish an \textit{a priori} bound for positive solutions of the regularized problem $(P_{\mathcal{B},\varepsilon})$ in $\mathbb{R}\times V$ (Corollary \ref{cor:abouD}). 

We start with an \textit{a priori} bound on $\lambda\in\mathbb{R}$, uniformly in $\varepsilon \in [0,1]$. 
\begin{prop} \label{prop:bpara} 
Assume %\mpfnsb{For tech.comment 10}
\eqref{H_1g}, \eqref{H_2g}, \eqref{H_3} and \eqref{ab:posi}. Then, there exist $\overline{\lambda}, \overline{\varepsilon}>0$ 
 such that if $0\leq \varepsilon\leq \overline{\varepsilon}$, then 
$(P_{\mathcal{B},\varepsilon})$ has no positive solution for $|\lambda| \geq\overline{\lambda}$.\end{prop}

\noindent\textit{Proof.} 
Let us suppose first that we are not in the case \eqref{excepc}.
%\marginpar{modified}\mpfns{agree}.
Let $0\leq \varepsilon \leq 1$, and assume that $(P_{\mathcal{B},\varepsilon})$ has a positive solution $u$ for some $\lambda > 0$. Let $B$ be given by 
\eqref{ab:posi}, $\lambda_B > 0$ be the first eigenvalue of the problem
\begin{align} \label{Depr}
\begin{cases}
-\Delta \phi = \lambda \phi & \mbox{ in } B, \\
\phi = 0 & \mbox{ on } \partial B,
\end{cases}
\end{align}
and  $\phi_B \in C^2(\overline{B})$ be a positive eigenfunction associated to $\lambda_B$. We extend $\phi_B$ to $\overline{\Omega}$ by setting $\phi_B = 0$ in $\overline{\Omega}\setminus \overline{B}$, so that $\phi_B \in H^1_0(\Omega)$.
Since $u>0$ on $\overline{B}$ and $\frac{\partial \phi_B}{\partial \nu}< 0$ on $\partial B$, the divergence theorem yields that 
\[
\lambda_B \int_B \phi_B u =
\int_B -\Delta \phi_B u = \int_B \nabla \phi_B \nabla u - \int_{\partial B} \frac{\partial \phi_B}{\partial \nu} u >\int_B \nabla \phi_B \nabla u, 
\]
where $\nu$ is the outward unit normal to $\partial B$. 
On the other hand, we see that
\begin{align*}
\int_B \nabla u \nabla \phi_B = \int_B b(x)g(u)\phi_B + \lambda \int_B
a(x)(u+\varepsilon)^{q-1} F(u) \phi_B,
\end{align*}
where $q$ is given by \eqref{H_3}. It follows that
\begin{align} \label{phB}
\int_B u^q \phi_B \left\{ b(x)\frac{g(u)}{u^q} + \lambda a(x) \left( \frac{u}{u+\varepsilon}\right)^{1-q} \frac{F(u)}{u} - \lambda_B u^{1-q} \right\} < 0.
\end{align}
Now, for $(x,s)\in B \times (0, \infty)$, we set 
\begin{align*}
h(x,s) := b(x)\frac{g(s)}{s^q} + \lambda a(x) \left( \frac{s}{s+\varepsilon}\right)^{1-q}
\frac{F(s)}{s} - \lambda_B s^{1-q}.
\end{align*}
By \eqref{H_2g}, there exists $s_0 > 0$ such that 
\begin{align*}
g(s)\geq \frac{\lambda_B}{b_0} s \quad \mbox{ for } s> s_0, 
\end{align*}
where $b_0$ is from \eqref{ab:posi}. 
Hence, since $f(s)>0$ for $s>0$ and $a(x)\geq 0$ a.e. in $B$, we deduce that if $\lambda > 0$, $x \in B$
and $s > s_0$, then
\begin{align} \label{h:posi1}
h(x,s) \geq b(x)\frac{g(s)}{s^q} -  \lambda_B s^{1-q} 
\geq s^{1-q} \left( b_0\frac{g(s)}{s} - \lambda_B \right) \geq 0.
\end{align}
Let us now consider the case $0<s\leq s_0$. From \eqref{H_1g}, we can choose $K_0 > 0$ such that
\[
\left| \frac{g(s)}{s} \right| \leq K_0 \quad\mbox{ for } 0<s\leq s_0.
\]
Recalling \eqref{H_3} (or \eqref{H_3'}), we set
\begin{align*}
M_1 := \inf_{0<s\leq s_0} \frac{F(s)}{s} > 0.
\end{align*}
By putting $b_\infty := \| b \|_\infty$, it follows that
\begin{align*}
h(x,s) & = s^{1-q} \left\{
b(x) \frac{g(s)}{s} + \lambda a(x) \left( \frac{1}{s+\varepsilon} \right)^{1-q}\frac{F(s)}{s} - \lambda_B \right\} \\
&\geq s^{1-q} \left\{ \lambda a_0 \left( \frac{1}{s_0+ 1}\right)^{1-q} M_1 - (\lambda_B + b_\infty K_0) \right\},
\end{align*}
so that $h(x,s) \geq 0$ for $x \in B$ and $0<s\leq s_0$ if
\begin{align} \label{h:posi2}
\lambda \geq \overline{\lambda}:=\frac{(\lambda_B + b_\infty K_0) (s_0 + 1)^{1-q}}{a_0M_1}.
\end{align}
Consequently, by \eqref{phB}, \eqref{h:posi1} and \eqref{h:posi2} we 
deduce that $\lambda < \overline{\lambda}$. 

Next, let us verify the existence 
of a lower bound on $\lambda < 0$  for the existence of a positive solution of $(P_{\mathcal{B},\varepsilon})$. In order to check this, we notice that if 
$(P_{\mathcal{B},\varepsilon})$ has a positive solution $u$ for some $\lambda < 0$ and $\varepsilon \in [0,1]$, then
\begin{align*}
-\Delta u = (-\lambda)(-a(x))(u+\varepsilon)^{q-1}F(u) + b(x) g(u) \quad\mbox{in $\Omega$}.
\end{align*}
From \eqref{ab:posi}, the desired conclusion follows 
arguing as above with $B$ now replaced by $B'$.

It remains to consider case \eqref{excepc}. However, it suffices to note that \eqref{ab:posi} implies that if $\varepsilon$ is small enough, 
then 
\begin{align*}
a_{\varepsilon}(x) \geq \frac{a_0}{2} \ \mbox{ in } B, \quad 
\mbox{ and } 
-a_{\varepsilon}(x) \geq a_0' \ \mbox{ in } B'. 
\end{align*}
The proof now follows in the same way as above. \qed \newline

Next, given a compact interval $I$, we establish an \textit{a priori} upper bound for positive solutions of $(P_{\mathcal{B},\varepsilon})$ whenever $\lambda \in I$ and $\varepsilon \in [0,1]$. We start with the following preliminary lemma
%\marginpar{modified}\mpfns{agree}
(see also \cite[Theorem 4.1]{ALG98}).

\begin{lem}  \label{lem:bound:norm}
Assume %\mpfnsb{Although in tech.comment 10, we remain the current approach using ABC.}
\eqref{H_2g}, \eqref{H_3}, \eqref{f:cave}, \eqref{g:vex} and
$(\mathcal{H}_b)$. Let $\Lambda > 0$. Suppose there exists a constant $C_1 > 0$ such that $\Vert u \Vert_{C(\overline{\Omega^b_+})} \leq C_1$ for all positive solutions $u$ of $(P_{\mathcal{B},\varepsilon})$ with $\lambda \in [-\Lambda, \Lambda]$ and $\varepsilon \in [0, 1]$. Then, there exists $C_2>0$ such that $\Vert u \Vert_{C(\overline{\Omega})} \leq C_2$ for all positive solutions $u$ of $(P_{\mathcal{B},\varepsilon})$ with $\lambda \in [-\Lambda, \Lambda]$ and $\varepsilon \in [0, 1]$.  
\end{lem}

\noindent\textit{Proof.}
(i) First, we consider the Dirichlet case. We use a comparison principle for concave problems inspired by the one in \cite[Lemma 3.3]{ABC94}. 

(1) Assume first that $\overline{\Omega^b_+} \subset \Omega$, 
and recall that $D_b$ is given by $(\mathcal{H}_b)$.   
Let $\lambda \in [0, \Lambda]$, and consider the problem
\begin{align} \label{p:a-}
\begin{cases}
-\Delta v = \lambda a^+ (x) (v+\varepsilon)^{q-1}F(v) - b^-(x) g(v) 
& \mbox{ in } D_b, \\ 
v =  C_1 & \mbox{ on } \partial \Omega^b_+, \\
v = 0 & \mbox{ on } \partial \Omega.  
\end{cases}
\end{align}
Let $u$ be a positive solution of $(P_{\mathcal{B},\varepsilon})$, with $\lambda\in [0, \Lambda]$ and $\varepsilon \in [0, 1]$. It follows  that $u>0$ in $D_b$. Since $\lambda \geq 0$, $b=-b^-$ in $D_b$, and $f(s)>0$ for $s>0$, 
it is easy to check that $u$ is a subsolution of \eqref{p:a-}.   

Next, we construct a supersolution of  \eqref{p:a-}. Consider
% \marginpar{replaced $w$ by $w_0$ below}\mpfns{I think $w$ is good.}
the unique positive solution $w_0$ of the problem 
$$
\begin{cases}
-\Delta w = 1 & \mbox{in $D_b$}, \\
w=0 & \mbox{on $\partial \Omega^b_+ \cup \partial \Omega$}. 
\end{cases}
$$
Set $\overline{w} := C(w_0+1)$ for $C>0$. If $C\geq C_1$, then 
$\overline{w} \geq C_1$ on $\partial \Omega^b_+$, and 
$\overline{w} \geq 0$ on $\partial \Omega$. Moreover, we claim that if $C$ is sufficiently large, then 
\begin{align*}
-\Delta \overline{w} \geq \lambda a^+(x) (\overline{w} 
+ \varepsilon)^{q-1} F(\overline{w})  -b^-(x) g(\overline{w}) 
\quad\mbox{ in } D_b. 
\end{align*}
Indeed, let 
\[
\delta := \frac{1}{\Lambda \| a^+ \|_\infty (\| w_0 \|_{C(\overline{D_b})}+1)} > 0. 
\]
Then, from \eqref{f:cave} and \eqref{H_2g} there exists $s_1 > 0$ large 
enough such that if $s\geq s_1$, then 
\begin{align*}
& 0 \leq g(s),  \mbox{ and } f(s) \leq \delta s \ (\mbox{so that } 
F(s) \leq \delta s^{2-q}).  
\end{align*}
It follows that if $C\geq s_1$, then 
%\mpfnsb{modified below}
\begin{align*}
&-\Delta \overline{w} - \left\{ \lambda a^+(x) (\overline{w} 
+ \varepsilon)^{q-1} F(\overline{w})  -b^-(x) g(\overline{w}) \right\}  \\ 
& \geq C - \Lambda \| a^+ \|_\infty \overline{w}^{q-1} \delta \overline{w}^{2-q}\\ 
& \geq C \left\{ 1 - \delta \Lambda \| a^+ \|_\infty \left( \| w_0 \|_{C(\overline{D_b})} + 1 \right) \right\} = 0 \quad \mbox{ in } D_b. 
\end{align*}
Thus the claim has been verified, and $\overline{w}$ is a supersolution of  \eqref{p:a-} if $C\geq \max (C_1, s_1)$. Note that $C$ can be chosen independently of 
$\lambda\in [0,\Lambda]$ and $\varepsilon \in [0, 1]$.

Now, we see from  
\eqref{f:cave} and \eqref{g:vex} that the nonlinearity in \eqref{p:a-} 
is concave, that is, if we set 
\begin{align*}
j(x, s) := \lambda a^+(x) (s + \varepsilon)^{q-1} F(s) - b^{-}(x) g(s), \quad 
x \in D_b, \ s > 0, 
\end{align*}
then 
\begin{align*}
\frac{\partial}{\partial s}\left( \frac{j(x,s)}{s} \right) < 0, 
\quad \mbox{for } x \in D_b, \ s>0. 
\end{align*}
Indeed,%\mpfnsb{modified}
\begin{align*}
\frac{d}{ds}\left( \frac{(s+\varepsilon)^{q-1}F(s)}{s} \right) 
&= (q-1)(s+\varepsilon)^{q-2}(s^{-q}f(s)) + (s+\varepsilon)^{q-1} \frac{d}{ds} 
(s^{-q}f(s)) \\
& < 0. 
\end{align*}
Reasoning as in \cite[Proposition A.1]{RQUTMNA} (whose argument is based on \cite[Lemma 3.3]{ABC94}), we may deduce that 
$u\leq \overline{w}$ in $D_b$, so that 
\[
u\leq C_1 + \| \overline{w} \|_{C(\overline{D_b})} \ \mbox{ on } \overline{\Omega}. 
\]

It remains to verify the case $-\Lambda \leq \lambda < 0$. Note that any  positive solution $u$ of $(P_{\mathcal{B},\varepsilon})$ with $\lambda \in [-\Lambda, 0)$ and $\varepsilon \in [0, 1]$ satisfies 
\begin{align*}
-\Delta u =  \mu (-a(x)) (u+\varepsilon)^{q-1}F(u) + b(x) g(u) \quad \mbox{in $D_b$},
\end{align*}
with $\mu :=- \lambda \in (0, \Lambda]$. Instead of \eqref{p:a-}, we consider the following concave problem:
\begin{align} \label{prb:mu-}
\begin{cases}
-\Delta v = \mu a^-(x) (v+\varepsilon)^{q-1}F(v)  - b^-(x) g(v) 
& \mbox{ in } D_b, \\
v =  C_1 & \mbox{ on } \partial \Omega^b_+, \\
v = 0 & \mbox{ on } \partial \Omega. 
\end{cases}
\end{align}
Then, we see that $u$ is a subsolution of this problem. 
The remainder of the argument is identical to the one in the case $\lambda \in [0, \Lambda]$. 

(2) Assume now that 
$\Omega^b_+ \supset \{ x \in \Omega : d(x, \partial \Omega) < \sigma \}$ for some $\sigma > 0$. This case can be verified identically. Indeed, it suffices to replace \eqref{p:a-} by the problem 
\begin{align*} 
\begin{cases}
-\Delta v = \lambda a^+ (x) (v+\varepsilon)^{q-1}F(v) - b^-(x) g(v) 
& \mbox{ in } D_b, \\ 
v =  C_1 & \mbox{ on } \partial \Omega^b_+ \cap \Omega. 
\end{cases}
\end{align*}

(ii) Lastly, we verify the Neumann case.
% \marginpar{removed "indeed"}\mpfns{agree}
Since $a_\varepsilon \leq a^+$ and $-a_\varepsilon \leq a^- + 1$, it suffices to replace \eqref{p:a-} by 
\begin{align*} 
\begin{cases}
-\Delta v = \lambda a^+ (x) (v+\varepsilon)^{q-1}F(v) - b^-(x) g(v) 
& \mbox{ in } D_b, \\ 
v =  C_1 & \mbox{ on } \partial \Omega^b_+, \\
\frac{\partial v}{\partial \mathbf{n}} = 0 & \mbox{ on } \partial \Omega, 
\end{cases}
\end{align*}
and \eqref{prb:mu-} by 
\begin{align*} 
\begin{cases}
-\Delta v = \mu (a^-(x) + 1) (v+\varepsilon)^{q-1}F(v)  - b^-(x) g(v) 
& \mbox{ in } D_b, \\
v =  C_1 & \mbox{ on } \partial \Omega^b_+, \\
\frac{\partial v}{\partial \mathbf{n}} = 0 & \mbox{ on } \partial \Omega. 
\end{cases}
\end{align*} 

The proof of Lemma \ref{lem:bound:norm} is now complete. \qed \newline

The following result is due to Gidas and Spruck \cite[Theorem 1.1]{GS81} and 
Amann and L\'opez-G\'omez \cite[Section 4]{ALG98} (see also L\'{o}pez-G\'{o}mez, Molina-Meyer and Tellini \cite[Section 6]{LGMMT13}).

\begin{prop}
\label{prop:bnorD} 
Assume %\mpfnsb{For tech.comment 10}
\eqref{H_2f}, \eqref{H_3} and \eqref{H_2'}. 
In addition, suppose either
\begin{enumerate}
\item[(i)] $b >0$ on $\overline{\Omega}$, and $p< \frac{N+2}{N-2}$ for the Dirichlet case ($p<\frac{N+1}{N-1}$ for the Neumann case) if $N>2$, or 
\item[(ii)] \eqref{f:cave}, \eqref{g:vex}, and $(\mathcal{H}_b)$ hold. 
\end{enumerate}
Then, given $\Lambda> 0$, there exists $C>0$ such that if $0\leq\varepsilon\leq1$, then $u\leq C$ in $\Omega$ for all positive  solutions $u$ of $(P_{\mathcal{B},\varepsilon})$ with $|\lambda|\leq\Lambda$.
\end{prop}

\begin{rem} 
\strut 
\begin{enumerate}
\item The condition (i) for the Neumann case is based on $(\mathcal{H}_b)$ with $\gamma = 0$ from Amann and L\'opez-G\'omez \cite[Section 4]{ALG98}, not on Gidas and Spruck \cite[Theorem 1.1]{GS81}. 

\item Under the condition (ii), we can handle the case where $b>0$ in $\Omega$, and $b=0$ somewhere on $\partial \Omega$. 
\end{enumerate}
\end{rem}

\noindent\textit{Proof.}
Case (i) under $\mathcal{B}u=u$ is verified in the same way as in the proof of \cite[Theorem 1.1]{GS81}.  

Case (i) under $\mathcal{B}u= \frac{\partial u}{\partial \mathbf{n}}$ and 
case (ii) can be handled using the arguments in \cite[Proposition 6.5]{RQUIsrael} and \cite[Proposition 4.2]{RQUTMNA}, which are based on \cite[Section 6]{LGMMT13}. Indeed, based on \eqref{H_2f}, \eqref{H_2'}, and $(\mathcal{H}_b)$,
we employ the argument developed by Amann and L\'opez-G\'omez \cite{ALG98} 
to deduce that the hypothesis of Lemma \ref{lem:bound:norm} is fulfilled, 
so that the desired conclusion follows. \qed \newline

By $\| \cdot \|_H$ we denote the usual norm of $H=H^1_0(\Omega)$ when $\mathcal{B}u=u$ and of $H=H^1(\Omega)$ when $\mathcal{B}u=\frac{\partial u}{\partial \mathbf{n}}$. The next %\mpfnsb{For tech.comment 11}
lemma follows easily by a bootstrap argument based on elliptic regularity and the Sobolev embedding theorem: 

\begin{lem}
\label{lem:norm} Assume that there exist $C>0$ 
and $p \in\left(  1, \frac{N+2}{N-2} \right)  $ such that
\begin{align*}
|g(s)| \leq C (1 + s^{p}) \quad\mbox{ for } s\geq0.
\end{align*}
Let $\Lambda > 0$ and $u$ be a nonnegative solution of $(P_{\mathcal{B},\varepsilon})$ for $\lambda \in [-\Lambda, \Lambda]$ and $\varepsilon \in 
[0,1]$. Then, given $c_1 > 0$ there exists $c_2 > 0$ such that 
if $\| u \|_{H} \leq c_1$, then $\| u\|_{V} \leq c_2$. 
\end{lem}

With the aid of Lemma \ref{lem:norm}, Propositions \ref{prop:bpara} and
\ref{prop:bnorD} provide us with an \textit{a priori} bound in 
$\mathbb{R}\times V$ for positive solutions of $(P_{\mathcal{B},\varepsilon})$ uniformly in $\varepsilon\in [0,1]$.

\begin{cor}
\label{cor:abouD} 
Under %\mpfnsb{For tech.comment 10}
the assumptions of Propositions \ref{prop:bpara} and \ref{prop:bnorD}, there exists $\overline{\varepsilon}, \overline{C}>0$ such that if $0\leq\varepsilon\leq \overline{\varepsilon}$, then $|\lambda|+\| u \|_{V} \leq\overline{C}$ for all positive solutions $u$ of $(P_{\mathcal{B},\varepsilon})$. 
\end{cor}

\noindent\textit{Proof.}
Propositions \ref{prop:bpara} and \ref{prop:bnorD} imply that there exist $\overline{\varepsilon}, C > 0$ such that if $0\leq \varepsilon \leq \overline{\varepsilon}$, then $|\lambda| + \| u \|_{H} \leq C$ for all positive solutions $u$ of  $(P_{\mathcal{B},\varepsilon})$, since $u$ satisfies
\[
\int_\Omega |\nabla u|^2 = \lambda \int_\Omega \{
a(x) (u+\varepsilon)^{q-1}F(u)u + b(x) g(u)u \} \leq C'
\]
for some $C'>0$. Lemma \ref{lem:norm} provides then the desired conclusion. \qed \newline

We discuss now bifurcation of nontrivial nonnegative solutions for $(P_{\mathcal{B}})$ from $(\lambda,0)$ for $\lambda > 0$. We prove the following preliminary lemma. 

\begin{lem}
\label{lem:boublw}
Assume \eqref{H_3} and $(H_\psi)$ with $\psi = a$. 
Let $\Lambda>0$, $\Omega^{\prime}$ be a nonempty connected open subset of
$\Omega^{a}_{+}$, and $B\Subset\Omega^{\prime}$ be a ball. Then there exists
$C_{\Lambda}>0$ such that if $\lambda\geq\Lambda$, then $\Vert u\Vert
_{C(\overline{B})}\geq C_{\Lambda}$ for all nontrivial nonnegative solutions
$u$ of $(P_{\mathcal{B}})$ such that $u\not \equiv 0$ in $\Omega^{\prime}$.
\end{lem}

%\begin{rem}
%	Indeed, 
%	the assumption \eqref{H_3} can be weakened by \eqref{f/sinfty} and \eqref{f:slope} (Remark \ref{rem:H3}).  
%\end{rem}

\noindent\textit{Proof.}
We use an argument based on sub and supersolutions. 
First of all, we remark that \eqref{H_3} implies \eqref{f/sinfty} and \eqref{f:slope}, see Remark \ref{rem:H3}. 

Let $\Lambda, \Omega'$ and $B$ be as in the statement of this lemma, and let $\lambda \geq \Lambda$. Assume that
$u$ is a nontrivial nonnegative solution of $(P_{\mathcal{B}})$ such that $u\not\equiv 0$ in $\Omega'$. Set 
\begin{align} \label{K2}
K_1 := \max_{0\leq s \leq \| u\|_{C(\overline{\Omega})}} |g'(s)| \geq 0.
\end{align}
Since $g(0)=0$, the mean value theorem provides some constant $\theta=\theta_x \in (0,1)$ such that $g(u)=g'(\theta u)u$. Thus, using \eqref{K2} we get that
\begin{align*}
& (-\Delta + b_\infty K_1 + 1)u \\
& = \lambda a(x)f(u) + (b_\infty K_1 + 1 + b(x) g'(\theta u))u \geq u \geq 0 \ \mbox{ and }
\not\equiv 0 \quad \mbox{ in } \Omega'.
\end{align*}
The strong maximum principle yields that $u>0$ in $\Omega'$.
Now, let $s_0 > 0$ be fixed. Then, the following two possibilities may occur: (i) $u \leq s_0$ on $\overline{B}$; (ii) $u > s_0$ somewhere on $\overline{B}$. 

We consider case (i). We have $a_0 := \inf_{\overline{B}} a > 0$, so that $u$ is a supersolution of the problem
\begin{align} \label{apr:cave}
\begin{cases}
-\Delta v + b_\infty K_{2} v = \lambda a_0 f(v) & \mbox{ in } B, \\
v = 0 & \mbox{ on } \partial B,
\end{cases}
\end{align}
where
\begin{align*} 
K_{2} := \max_{0\leq s \leq s_0} |g'(s)| \geq 0.
\end{align*}
Indeed, $u\geq 0$ on $\partial B$. Moreover, since $f(s)>0$ for $s>0$, the mean value theorem shows that
\begin{align*}
-\Delta u + b_\infty K_{2} u - \lambda a_0 f(u)
& = \lambda (a(x) - a_0) f(u) + (b_\infty K_{2} + b(x) g'(\theta u))u \\
& \geq 0 \quad \mbox{ in } B.
\end{align*}
To construct a subsolution of \eqref{apr:cave}, we use the positive eigenfunction $\phi_B$ associated to the first eigenvalue $\lambda_B$ of \eqref{Depr} and such that
$\| \phi_B \|_{C(\overline{B})} = 1$. From \eqref{f/sinfty}, we find a constant $s_1 > 0$ small enough such that
\begin{align} \label{Z2}
\frac{f(s)}{s} \geq \frac{\lambda_B + b_\infty K_2}{\Lambda a_0} \quad \mbox{ for } 0 < s \leq s_1.
\end{align}
If $0<s\leq s_1$, then we observe that 
\begin{align*}
-\Delta (s \phi_B) + b_\infty K_{2} s \phi_B - \lambda a_0 f(s \phi_B)
& \leq  s \phi_B \left\{ \lambda_B + b_\infty K_2 - \Lambda a_0
\frac{f(s \phi_B)}{s \phi_B} \right\} \\
& \leq 0 \quad \mbox{ in } B.
\end{align*}
This implies that $s \phi_B$ is a subsolution of \eqref{apr:cave} whenever
$0< s \leq s_1$.
Now, since $u > 0$ in $\Omega'$, it follows that $u>0$ on $\overline{B}$. Furthermore, we assert that
\begin{align} \label{ucLam}
u \geq s_1 \phi_B \quad\mbox{ on } \overline{B}.
\end{align}
By contradiction, we assume that $u \not\geq s_1 \phi_B$. 
Then, since $u>0 = s_1
\phi_B$ on $\partial B$, we can choose $\sigma \in (0,1)$ such that $u - \sigma
s_1 \phi_B \geq 0$ on $\overline{B}$ and $u - \sigma s_1 \phi_B = 0$ somewhere in $B$. From \eqref{f:slope}, we fix $M_{0} > 0$ such that 
\begin{align*}
\frac{f(s) - f(t)}{s-t} > - M_{0} \quad\mbox{ for } 0\leq t < s\leq s_{0}.
\end{align*}
Putting $M_1 := \lambda a_0 M_0 > 0$, we see that the mapping
\[
s \longmapsto M_1 s + \lambda a_0 f(s)
\]
is nondecreasing for $0\leq s \leq s_0$. Indeed, if $0\leq t < s \leq s_0$, then
\begin{align*}
& M_1 s + \lambda a_0 f(s) - (M_1 t + \lambda a_0 f(t)) \\
&= \left( M_1 + \lambda a_0 \frac{f(s) - f(t)}{s-t} \right)(s-t) \\
& \geq \left( M_1 - \lambda a_0 M_0 \right)(s-t) = 0,
\end{align*}
as desired.
Thus, using this monotonicity and having in mind that 
\[
(-\Delta + b_\infty K_2) \sigma s_1 \phi_B \leq \lambda a_0 f(\sigma s_1 \phi_B) 
\ \mbox{ in } B \quad (\mbox{recall \eqref{Z2}}), 
\]
we deduce that
\begin{align*}
& (-\Delta + b_\infty K_2 + M_1) (u - \sigma s_1 \phi_B) \\
&\geq M_1 u + \lambda a_0 f(u) - \left( M_1 \sigma s_1 \phi_B + \lambda
a_0 f(\sigma s_1 \phi_B) \right) \\
&\geq 0 \quad \mbox{ in } B,
\end{align*}
and $u - \sigma s_1 \phi_B = u > 0$ on $\partial B$. The strong maximum 
principle yields that $u - \sigma s_1 \phi_B > 0$ in $B$, a contradiction. Thus, we have verified  \eqref{ucLam}. By taking into account case (ii), $C_{\Lambda}: = \min \{ s_0, s_1 \}$ is as desired. \qed \newline

By virtue of Lemma \ref{lem:boublw}, there are no nontrivial nonnegative solutions of $(P_{\mathcal{B}})$ bifurcating from $(\lambda,0)$ for $\lambda>0$, and moreover, there exist no small positive solutions of $(P_{\mathcal{B}})$ for $\lambda=0$. In view of this fact, 
although we shall observe that $(P_{\mathcal{B}})$ possesses 
a bounded subcontinuum of nontrivial nonnegative solutions bifurcating at
$(0,0)$, we infer that the bifurcation subcontinuum 
is of \textit{loop type}.

\begin{prop}
\label{prop:nobifz}
Assume \eqref{H_1g}. %\mpfnsb{For tech.comment 12}
Then, the following assertions hold:

\begin{enumerate}
\item Assume %\mpfns{modified}
$(H_\psi)$ with $\psi = a$, and additionally \eqref{H_3} if $\mathcal{B}u=\frac{\partial u}{\partial \mathbf{n}}$. 
Then, given $\lambda_0> 0$, there exist $\delta_{0},
c_{0}> 0$ such that $\| u \|_{C(\overline{\Omega})} \geq c_{0}$
for all nontrivial nonnegative solutions $u$ of $(P_{\mathcal{B}})$ for 
$\lambda\in(\lambda_{0}- \delta_{0}, \lambda_{0}+ \delta_{0})$. 

\item Assume %\mpfns{modified}
\eqref{ib<0} and \eqref{gssig} if $\mathcal{B}u=\frac{\partial u}{\partial \mathbf{n}}$. Then, there exists $C>0$ such that $\|
u \|_{C(\overline{\Omega})}\geq C$ for all positive solutions $u$ of 
$(P_{\mathcal{B}})$ for $\lambda=0$. 
\end{enumerate} 
\end{prop}

%	\begin{rem}
%	The \mpfns{In the same reason we put this remark.} 
%	assumption \eqref{H_3} can be weakened by \eqref{f/sinfty} and 
%	\eqref{f:slope} if $\mathcal{B}u=u$, 
%	whereas it can be weakened by \eqref{H_3'} and \eqref{f:slope} if 
%	$\mathcal{B}u=\frac{\partial u}{\partial \mathbf{n}}$ (Remark \ref{rem:H3}).
%	\end{rem}

%\mpfns{Remark deleted. Instead, add the remark at the beginning of the proof.}

\noindent\textit{Proof.}
(i) %\marginpar{modified}\mpfns{agree}
We recall that \eqref{H_3} also implies \eqref{H_3'}, see Remark \ref{rem:H3}. First, we verify the Dirichlet case. 
By contradiction, we assume that $\lambda_n \to \lambda_0 > 0$, and $u_n$ are nontrivial nonnegative solutions of $(P_{\mathcal{B}})$ with $\lambda=\lambda_n$ such that
$\| u_n \|_{C(\overline{\Omega})} \to 0$. Then, we claim that, up to a subsequence,
$\int_\Omega a(x) f(u_n)u_n \leq 0$. If not, then we may suppose that $\int_\Omega a(x) f(u_n) u_n > 0$ for all $n$. It follows that $u_n \not\equiv 0$ in 
$\Omega^a_+$. Indeed, if $u_n\equiv 0$ in $\Omega^a_+$, then, using that $f(s)>0$ for $s>0$, we find that
\begin{align*}
0< \int_\Omega a(x) f(u_n) u_n \leq \int_{\Omega^a_+} a(x) f(u_n) u_n = 0,
\end{align*}
which is a contradiction. 
Employing $(H_\psi)$ with $\psi=a$, we may deduce that there exists a connected open subset $\Omega' \subset \Omega^a_+$ such that
$u_n \not\equiv 0$ in $\Omega'$ for all $n\geq 1$.
Let $B \Subset \Omega'$ be a ball. 
We apply Lemma \ref{lem:boublw} with
$\Lambda = \frac{\lambda_0}{2}$, to derive that
$\| u_n \|_{C(\overline{B})} \geq c_0$ for some $c_0 > 0$ 
independent of $n$, which contradicts $\| u_n \|_{C(\overline{\Omega})} 
\to 0$. Thus, the claim follows.

Now, we observe from the definition of $u_n$ that
\begin{align*}
\| u_n \|_{H}^2 :=
\int_\Omega |\nabla u_n|^2  &= \lambda_n \int_\Omega a(x) f(u_n)u_n + \int_\Omega 
b(x)g(u_n)u_n \\ 
& \leq \int_\Omega b_\infty |g(u_n)| u_n.
\end{align*}
Set $v_n := u_n/ \| u_n \|_{H}$, so that
\begin{align} \label{vnXleq}
\| v_n \|_{H}^2 \leq b_\infty \int_\Omega \frac{|g(u_n)| u_n}{\| u_n \|_{H}^2}.
\end{align}
From \eqref{H_1g}, given $\varepsilon > 0$ there exists $s_\varepsilon > 0$ such that 
\[
|g(s)| \leq \frac{\varepsilon}{b_\infty} s \quad \mbox{ for } 0<s\leq s_\varepsilon.  
\]
Also, for $n$ large enough, we have that $\| u_n \|_{C(\overline{\Omega})} \leq s_\varepsilon$, so that
\begin{align} \label{leqepvn} 
b_\infty \int_\Omega \frac{|g(u_n)| u_n}{\| u_n \|_{H}^2} \leq \varepsilon \int_\Omega v_n^2.
\end{align}
From \eqref{vnXleq} and \eqref{leqepvn}, we derive that $v_n \to 0$ in $H^1_0(\Omega)$, a contradiction. \newline 

Next, we verify the Neumann case. Assume to the contrary that $\lambda_n \to \lambda_0 > 0$, and $u_n$ are nontrivial nonnegative solutions  of $(P_{\mathcal{B}})$ with $\lambda=\lambda_n$ such that $\| u_n \|_{C(\overline{\Omega})} 
\to 0$. 
We remark that $\| u_n \|_{H} \to 0$, since $u_n$ are nonnegative solutions of $(P_{\mathcal{B}})$ with $\lambda=\lambda_n$. Arguing as in the proof 
for the Dirichlet case, we have that, up to a subsequence, 
$\int_\Omega a(x) f(u_n)u_n \leq 0$, and consequently, 
$\int_\Omega |\nabla u_n|^2 \leq b_\infty \int_\Omega |g(u_n)| u_n$. 
%\mpfns{deleted}
%, where $b_\infty = \| b \|_\infty$. 
\par
Set $v_n := u_n/\| u_n \|_{H}$, so that $\| v_n \|_{H} = 1$.
We may assume that there exists $v_0 \in H^1(\Omega)$ such that $v_n \rightharpoonup v_0$ in $H^1(\Omega)$, $v_n \to v_0$ a.e. in $\Omega$, and $v_n \to v_0$ in $L^t(\Omega)$ for $t <2^*$. 
By \eqref{H_1g}, for any $\varepsilon > 0$ there exists $s_0 > 0$ such that
\begin{align} \label{g:ep}
|g(s)| \leq \frac{\varepsilon}{b_\infty} s \quad \mbox{ for } 0\leq s\leq s_0.
\end{align}
Thus, for $n$ large enough we have that $\| u_n \|_{C(\overline{\Omega})}
\leq s_0$, so that
\begin{align*}
\int_\Omega |\nabla v_n|^2 \leq b_\infty \int_\Omega \frac{|g(u_n)|}{\| u_n \|_{H}} v_n \leq \varepsilon \int_\Omega v_n^2 \leq \varepsilon.
\end{align*}
This implies that $\int_\Omega |\nabla v_n|^2 \to 0$, and it follows that $v_n \to v_0$, and $v_0$ is a positive constant.

Since $u_n$ is a nonnegative solution of $(P_{\mathcal{B}})$ with $\lambda=\lambda_n$, we see that, for every $\phi \in H^1(\Omega)$, 
\begin{align} \label{1-qfg}
\left( \int_\Omega \nabla v_n \nabla \phi \right) \| u_n \|_{H}^{1-q}
= \lambda_n \int_\Omega a(x) \frac{f(u_n)}{\| u_n \|_{H}^q} \phi 
+ \int_\Omega b(x)\frac{g(u_n)}{\| u_n \|_{H}^q} \phi.
\end{align}
Since $\| u_n \|_{C(\overline{\Omega})} \to 0$, \eqref{g:ep}
implies that $|g(u_n)|\leq u_n$ for $n$ large enough, so that
%\mpfnsb{modified below}
\begin{align*}
\left| \int_\Omega b(x)\frac{g(u_n)}{\| u_n \|_{H}^q} \phi \right|
&\leq b_\infty \int _\Omega \frac{u_n}{\| u_n \|_{H}^q} |\phi| \\
& \leq b_\infty \| u_n \|_{H}^{1-q} \left( \int_\Omega v_n^2 \right)^{\frac{1}{2}}
\left( \int_\Omega \phi^2 \right)^{\frac{1}{2}} \longrightarrow 0.
\end{align*}
We use this inequality to deduce from \eqref{1-qfg} that, passing to 
the limit as $n\to \infty$,
\begin{align} \label{afph/unq}
\int_\Omega a(x) \frac{f(u_n)}{\| u_n \|_{H}^q} \phi \longrightarrow 0.
\end{align}

On the other hand, since $f(0)=0$ we have that   
\[
\int_\Omega a(x) \frac{f(u_n)}{\| u_n \|_{H}^q} \phi=\int_{v_n>0} a(x) \frac{f(\|u_n\|_{H}v_n)}{(\| u_n \|_{H}v_n)^q} v_n^q\phi. 
\]
Thus, using \eqref{H_3'} and the fact that $u_n \to 0$ in $H^1(\Omega)$, $v_n \to v_0$ in $L^t(\Omega)$, $v_n \to v_0$ a.e. in $\Omega$ and $v_0$ is a positive constant,  the Lebesgue dominated convergence theorem yields that 
\[
\int_\Omega a(x) \frac{f(u_n)}{\| u_n \|_{H}^q} \phi \longrightarrow 
\frac{f_0}{q} v_0^q \int_\Omega a(x) \phi. 
\]
Therefore
\[
\int_\Omega a(x) \phi = 0.
\]
Since $\phi \in H^1(\Omega)$ is arbitrary, 
we find that $a\equiv 0$, which is a contradiction. \newline 

(ii) In the Dirichlet case, we argue as in the proof of assertion (i) to prove assertion (ii), by taking $\lambda_n = \lambda_0 = 0$ therein. 

Next, we verify the Neumann case. 
Assume to the contrary that there exist positive solutions $u_n$ of 
$(P_{\mathcal{B}})$ 
%\mpfnsb{modified on 10jun}
for $\lambda=0$ such that $\| u_n \|_{C(\overline{\Omega})} \to 0$. Then, 
as in the proof of assertion (i), 
we may deduce from \eqref{H_1g} that $v_n := u_n/ \| u_n\|_{H^1(\Omega)} \to v_0$ in $L^t (\Omega)$ for 
$t <2^*$, and $v_0$ is a positive constant. Since $u_n$ is a positive solution of $(P_{\mathcal{B}})$ %\mpfnsb{modified on 10jun}
with $\lambda=0$, we obtain $\int_\Omega b(x)g(u_n) = 0$. Recalling \eqref{gssig}, we see that 
\begin{align*}
0 = \int_\Omega b(x) \frac{g(u_n)}{\| u_n \|_{H^1(\Omega)}^\sigma} 
\longrightarrow g_0 v_0^\sigma \int_\Omega b(x), 
\end{align*}
so that $\int_\Omega b(x) = 0$, which contradicts \eqref{ib<0}.

The proof is now complete. \qed \newline

Assuming additionally $(H_\psi)$ with $\psi = -a$, 
we can extend %\mpfns{added}
Proposition \ref{prop:nobifz}(i) to $\lambda< 0$,
and in this case, bifurcation of nontrivial nonnegative solutions
of $(P_{\mathcal{B}})$ from $(\lambda,0)$ can only occur at $(0,0)$.

\begin{cor}
\label{cor:nobifzD} Under the assumptions of Proposition 
\ref{prop:nobifz}, assume in addition $(H_\psi)$ with $\psi = -a$. 
Then the conclusion of Proposition \ref{prop:nobifz}(i) 
holds for all $\lambda_0\neq0$. 
%\mpfnsb{For tech.comment 12} %\mpfns{modified}
In particular, given $\delta \in (0,1)$, the set of nontrivial nonnegative solutions of $(P_{\mathcal{B}})$ is away from the set $\{(\lambda,0); \delta \leq |\lambda|\leq \delta^{-1}\}$. 
\end{cor}

\noindent\textit{Proof.}
In view of Proposition \ref{prop:nobifz}, it remains to verify the case $\lambda < 0$. Assume to the contrary that $\lambda_n \to \lambda_0 < 0$, and $u_n$ are nontrivial nonnegative solutions of $(P_{\mathcal{B}})$ with $\lambda=\lambda_n$ such that $u_n \to 0$ in $C(\overline{\Omega})$. Then, we have that
\begin{align*}
-\Delta u_n = (-\lambda_n)(-a(x)) f(u_n) + b(x)g(u_n) \quad\mbox{in $\Omega$}.
\end{align*}
By the same arguments used in Lemma \ref{lem:boublw} and 
Proposition \ref{prop:nobifz}, we get the desired conclusion. \qed 

%----------------------------------------------------------
\section{Proofs of Theorems \ref{resu01} and \ref{resu02}} 
\label{sec:prf1} 
If $X$ is a 
%\mpfns{`complete' deleted} 
metric space and $E_n \subset X$, then we set
\begin{align*}
& \liminf_{n\to \infty} E_n := \{ x \in X : \lim_{n\to \infty}{\rm dist}\, (x, E_n) = 0 \}, \\
& \limsup_{n\to \infty} E_n := \{ x \in X : \liminf_{n\to \infty}{\rm dist}\, (x, E_n) = 0 \}.
\end{align*}
We shall use the following result due to Whyburn \cite[(9.12)Theorem]{Wh64}:

\begin{theorem}\label{thm:W}
Assume $\left\{
E_{n}\right\}$ is a sequence of connected sets satisfying that 
\begin{enumerate}
\item[(i)] $\displaystyle \bigcup_{n\geq 1}E_n$ is precompact; 
\item[(ii)] $\displaystyle \liminf_{n\to \infty}E_n \neq \emptyset$. 
\end{enumerate}
Then, $\displaystyle \limsup_{n\to \infty}E_n$ is nonempty, closed and 
connected. 
\end{theorem}

%	\begin{rem}
%	The original statement of Theorem \ref{thm:W} in \cite{Wh64} does not mention in its conclusion that 
%	$\limsup_{n\to \infty}E_n$ is closed. However, we can check the closedness in the following way. 
%
%	Assume $x_k \in \limsup_n E_n$ and $\mathrm{dist}(x_k, x_\infty) \to 0$. Then, we shall prove that $x_\infty \in \limsup_n E_n$. 
%	For $j\geq 1$, we choose a subsequence $\left\{
%	x_{k_j}\right\}$ such that $\mathrm{dist}(x_{k_j}, x_\infty) < \frac{1}{2j}$. Since $x_k \in \limsup_n E_n$, we also choose $y_{n_j} \in E_{n_j}$ such that $\mathrm{dist}(x_{k_j}, y_{n_j}) < \frac{1}{2j}$. It follows that 
%	\begin{align*}
%	\mathrm{dist}(x_\infty, E_{n_j}) &\leq \mathrm{dist}(x_\infty, y_{n_j}) \\ 
%	& \leq \mathrm{dist}(x_\infty, x_{k_j}) + \mathrm{dist}(x_{k_j}, y_{n_j}) 
%	< \frac{1}{2j} + \frac{1}{2j} = \frac{1}{j} \longrightarrow 0. 
%	\end{align*}
%	This implies that $\liminf_{n\to \infty}\mathrm{dist}(x_\infty, E_n) = 0$, as desired. \\ 
%	\end{rem}
%

As stated in Remark \ref{case>0}, 
we only have to prove Theorem \ref{resu02}(i) in the case
$\int_\Omega a < 0$. 
%	\begin{align} \label{ia>0}
%	\int_\Omega a < 0.
%	\end{align}  

\subsection{Proof of assertion (i) in Theorems \ref{resu01} and \ref{resu02}} 
\label{subsec:prfi}
When $\mathcal{B}u=\frac{\partial u}{\partial \mathbf{n}}$, we employ the following lemma, which concerns the direction of the bifurcation component $\mathcal{C}_\varepsilon$ at $(0,0)$, see \cite[Theorem 5.1]{RQUIsrael} for the proof.

\begin{lem} \label{prop:bifdir}
Let $\mathcal{B}u=\frac{\partial u}{\partial \mathbf{n}}$. 
Assume \eqref{H_3}, \eqref{gssig} and $\int_\Omega a\neq 0$. 
%\mpfns{Lemma \ref{prop:bifdir} modified}
% and $\int_\Omega a < 0$. %\eqref{ia>0}. 
Let $Z$ be any complement of $\langle 1 \rangle$ in $W^{2,r}(\Omega)$. Then, 
%\mpfns{modified}
for $\varepsilon > 0$ small enough, the set $\{ (\lambda, u) \}$ of nontrivial solutions of $(P_{\mathcal{B},\varepsilon})$ around $(0, 0)$ is parametrized as
$$
(\lambda, u) = ( \mu (s), \ s(1 + z (s))  ),
$$
with $s \in (-s_0, s_0)$, for some $s_0 >0$. Here $\mu: (-s_0, s_0) \to \mathbb{R}$ and $z: (-s_0, s_0) \to Z$ are continuous, and satisfy $\mu (0) =z (0)= 0$. Therefore, $\mathcal{C}_\varepsilon$ is precisely described by $\{ (\mu (s), s(1 + z (s))) : s \in [0, s_0) \}$ around $(0,0)$. Furthermore, the following holds:
\begin{align} \label{dir1348}
\lim_{s\to 0^+} \frac{\mu (s)}{s^{\sigma -1}} = - \varepsilon^{1-q} \frac{q g_0 \int_\Omega b(x)}{f_0 \int_\Omega a(x)}. 
\end{align}
In particular, under \eqref{ib<0} and $\int_\Omega a < 0$, 
the bifurcation of $\mathcal{C}_\varepsilon$ is {\rm subcritical} at $(0,0)$.
%$\mu (s)<0$ for $s>0$ small.
\end{lem}

Now, we consider the metric space $X := \mathbb{R}\times V$ 
with the metric given by %\mpfns{\textbf{For tech.comment 13}}
\[
d((\lambda, u), (\mu, v)) := |\lambda - \mu| + \| u - v \|_{V}\quad\mbox{for} \ \ (\lambda, u), (\mu, v) \in \mathbb{R}\times V.
\]
From Corollary \ref{cor:abouD}, if $\varepsilon \in (0,1]$, then
the components $\mathcal{C}_\varepsilon^\pm$ of positive solutions of $(P_{\mathcal{B},\varepsilon})$, emanating from
$(\lambda_{1,\varepsilon}^\pm, 0)$, satisfy
\begin{align} \label{Cep:bdd}
\mathcal{C}_\varepsilon^{\pm} \subset \{ (\lambda, u) \in \mathbb{R}\times V:|\lambda| + \| u \|_{V} \leq \overline{C} \}, 
\end{align}
where $\overline{C}$ does not depend on $\varepsilon \in (0,1]$. This implies that $\mathcal{C}_\varepsilon^\pm$ are both bounded, and consequently, we deduce that $\mathcal{C}_\varepsilon^- = \mathcal{C}_\varepsilon^+$ (cf. \cite[Proposition 18.1]{Am76}).  Then,
$\mathcal{C}_\varepsilon := \mathcal{C}_\varepsilon^\pm$ is nonempty and connected. In addition,
\begin{align} \label{liminfnon}
(0,0) \in \liminf_{\varepsilon \to 0^+} \mathcal{C}_\varepsilon,
\end{align}
since $\lambda_{1,\varepsilon}^\pm \to 0$ as $\varepsilon \to 0^+$. Moreover, by elliptic regularity, we obtain that
\begin{align} \label{precom}
\bigcup_{\varepsilon > 0} \mathcal{C}_{\varepsilon} \mbox{ is precompact. }
\end{align}
Indeed, for any $\{ (\lambda_n, u_n) \} \subset \bigcup_{\varepsilon > 0} \mathcal{C}_{\varepsilon}$ we have that $(\lambda_n, u_n) \in \mathcal{C}_{\varepsilon_n}$ for some $\varepsilon_n \in (0,1]$. From \eqref{Cep:bdd}, we may assume that $\{\lambda_n\}$ is a convergent sequence. Using \eqref{Cep:bdd} again, we deduce that $u_n \in W^{2,r}(\Omega)$ are solutions of
\begin{align*}
\begin{cases}
-\Delta u_n = \lambda_n a(x) (u_n + \varepsilon_n)^{q-1} F(u_n) + b(x)g(u_n) &
\mbox{ in } \Omega, \\
u_n = 0 & \mbox{ on } \partial \Omega.
\end{cases}
\end{align*}
In particular, using a bootstrap argument and the Sobolev embedding theorem, we deduce that  $\Vert u_n \Vert_{C^{1+\theta}(\overline{\Omega})}$ is bounded, for some $\theta \in (0,1)$.
The compact embedding $C^{1+\theta}(\overline{\Omega}) \subset C^{1}(\overline{\Omega})$ implies that $\{ u_n \}$ has a convergent subsequence in $V$, as desired.
Now, by \eqref{liminfnon} and \eqref{precom}, we may apply Theorem \ref{thm:W} to infer that $\mathcal{C}_{0} := \limsup_{\varepsilon \to 0^+} \mathcal{C}_{\varepsilon}$ is non-empty, closed and connected in $\mathbb{R}\times V$. 
From \eqref{Cep:bdd}, $\mathcal{C}_{0}$ is bounded in $\mathbb{R}\times V$. In addition, $\mathcal{C}_{0}$ is contained in the nonnegative solutions set of 
$(P_{\mathcal{B}})$. Indeed, given $(\lambda, u) \in \mathcal{C}_{0}$, there exists $(\lambda_n, u_n) \in \mathcal{C}_{\varepsilon_n}$ such that $\varepsilon_n \to 0^+$ and $(\lambda_n, u_n) \to (\lambda, u)$ in $\mathbb{R} \times V$. Thus $u$ is a nonnegative weak solution of $(P_{\mathcal{B}})$, and eventually, 
a nonnegative solution in $W^{2,r}(\Omega)$ by elliptic regularity.

Now, we show that $\mathcal{C}_{0}$ is nontrivial. By construction, we see that for $\varepsilon \to 0^+$, there exists a positive solution $u_\varepsilon$ of 
$(P_{\mathcal{B},\varepsilon})$ such that $(0, u_\varepsilon) \in \mathcal{C}_\varepsilon$. Indeed, 
%\mpfns{added}
we used Lemma \ref{prop:bifdir} if $\mathcal{B}u=\frac{\partial u}{\partial \mathbf{n}}$. In this case, we observe from \eqref{dir1348} that when 
%\mpfns{modified}%\eqref{ia>0} and 
\eqref{ib<0} and $\int_\Omega a < 0$ hold, the bounded component $\mathcal{C}_\varepsilon$ bifurcates %\mpfns{\textbf{For tech.comment 13}}
subcritically at $(0,0)$, provided that $\varepsilon$ is small enough. 
This implies that $\mathcal{C}_\varepsilon$ cuts $\{ (0, u) : 0\not\equiv u \geq 0 \}$, and consequently, the desired assertion follows. 
Since $\| u_\varepsilon \|_{V}\leq \overline{C}$, it follows by combining elliptic regularity and standard compactness arguments as above  that there exist $\varepsilon_n \to 0^+$  and $u_n := u_{\varepsilon_n}$ such that $u_n$ converges in $V$ to a nonnegative solution 
$u_0$ of $(P_{\mathcal{B}})$ for $\lambda=0$. By definition, we have that $(0, u_0) \in \mathcal{C}_{0}$. From Proposition \ref{prop:nobifz}(ii), 
we infer that $u_0$ is nontrivial, and so, $u_0 \gg 0$ by the 
strong maximum principle and Hopf's lemma. Assertion (i-1) has been 
now verified. We use Proposition \ref{prop:nobifz}(ii) again to deduce assertion (i-2).

Since $\mathcal{C}_0$ is nontrivial, we infer from Corollary \ref{cor:nobifzD}
that $\mathcal{C}_{0}$ does not contain any $(\lambda, 0)$ with $\lambda\neq 0$. Assertion \eqref{single} has been verified.

Finally, we verify assertion (i-3). 
For $\rho> 0$ and $(\lambda_{1}, u_{1}) \in\mathbb{R}\times V$, we set
\begin{align*}
&  B_{\rho}((\lambda_{1}, u_{1})) := \{ (\lambda, u) \in\mathbb{R} \times V:
|\lambda- \lambda_{1}| + \| u - u_{1}
\|_{V} < \rho\},\\ 
&  S_{\rho}((\lambda_{1}, u_{1})) := \{ (\lambda, u) \in\mathbb{R} \times V:
|\lambda- \lambda_{1}| + \| u - u_{1}
\|_{V} = \rho\}. 
\end{align*}
We note that $\overline{B_{\rho}((\lambda_{1}, u_{1}))} = B_{\rho}%
((\lambda_{1}, u_{1})) \cup S_{\rho}((\lambda_{1}, u_{1}))$.

\begin{figure}[!htb]
		\begin{center}
		\includegraphics[scale=0.16]{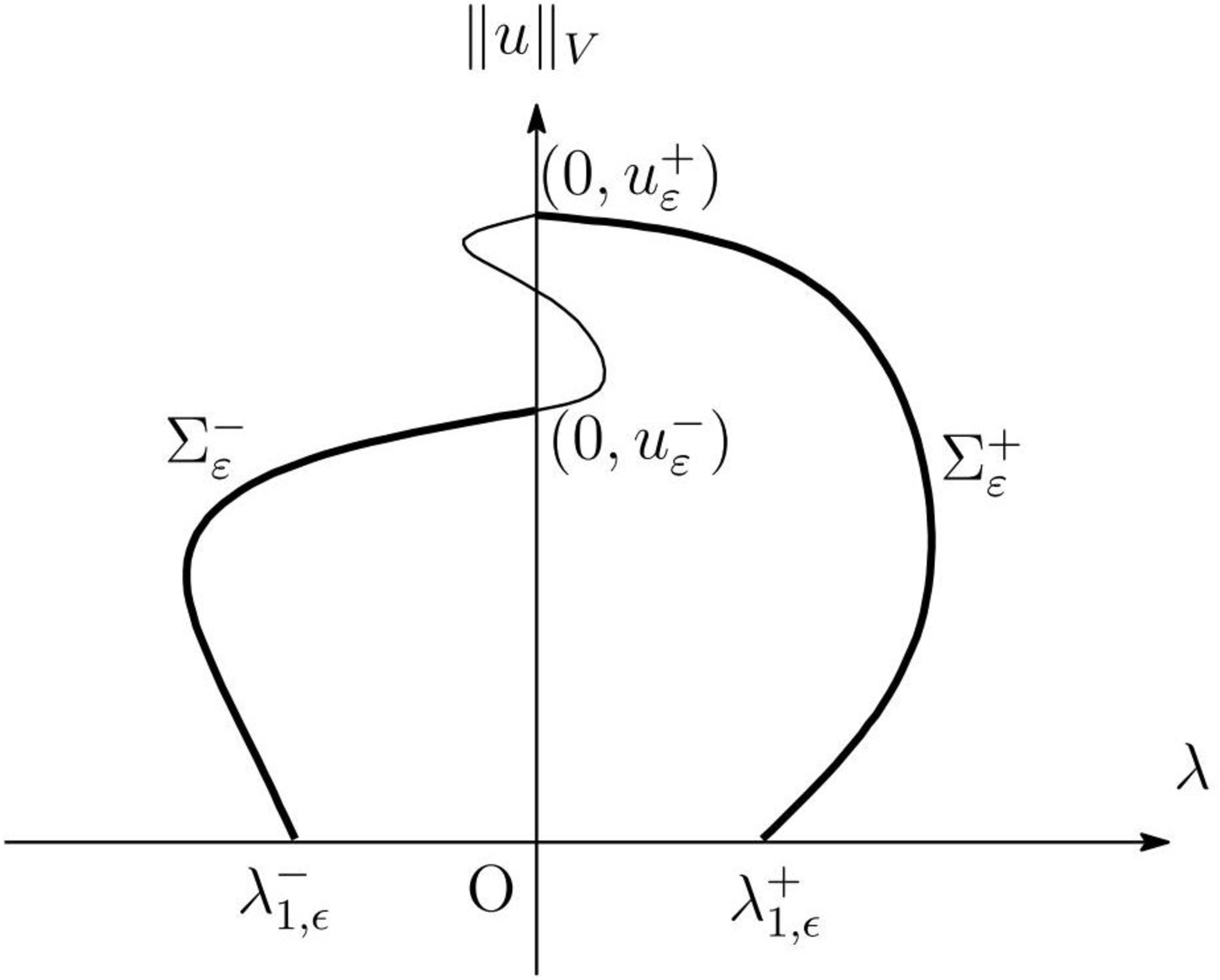} 
			\caption{The situations of 
				$\Sigma_{\varepsilon}^\pm$.} 
			\label{fig17_0302}
		\end{center}		  
		 \end{figure}

Let $\Sigma^{+}_{\varepsilon}$ and $\Sigma^{-}_{\varepsilon}$ be closed connected 
subsets of $\{ (\lambda, u) \in \mathcal{C}_\varepsilon : \lambda \geq 0 \}$ and 
$\{ (\lambda, u) \in \mathcal{C}_\varepsilon : \lambda \leq 0 \}$, respectively, 
such that $(\lambda_{1, \varepsilon}^{\pm}, 0), (0, u_\varepsilon^{\pm}) \in \Sigma^{\pm}_{\varepsilon}$ 
for some positive solutions $u_\varepsilon^{\pm}$ of 
$(P_{\mathcal{B}})$ for $\lambda=0$, see Figure \ref{fig17_0302}. 
This is well defined thanks to Proposition \ref{prop:nobifz}(ii). 
Since $\Sigma^\pm_{\varepsilon} \subset \mathcal{C}_\varepsilon$, we observe 
that
\begin{align} \label{C0rho}
\Sigma^\pm_{0} :=
\limsup_{\varepsilon \to 0^+} \Sigma^\pm_{\varepsilon} \subset
\limsup_{\varepsilon \to 0^+} \mathcal{C}_\varepsilon = \mathcal{C}_{0}.
\end{align}
Repeating the argument above, 
Whyburn's topological approach yields that
$\Sigma^\pm_{0}$ are non-empty, closed and connected sets consisting of nonnegative solutions of $(P_{\mathcal{B}})$ and such that $(0,0) \in \liminf_{\varepsilon \to 0^+} \Sigma^\pm_{\varepsilon} \subset \Sigma^\pm_{0}$. Proposition \ref{prop:nobifz}(ii) tells us that $(0, u_0^{\pm}) \in \Sigma_{0}^{\pm}$, for some positive solutions $u_0^{\pm}$ of $(P_{\mathcal{B}})$ with $\lambda=0$. 
It follows that  
$\Sigma^\pm_{0} \neq \{ (0,0) \}$, and by virtue of Corollary \ref{cor:nobifzD}, that 
$\Sigma^\pm_{0} \setminus \{ (0,0)\}$ consists of nontrivial nonnegative solutions of $(P_{\mathcal{B}})$. 

By definition, $(\lambda, u) \in \Sigma^+_{0}$ (respect.\ $\Sigma^-_{0}$) implies $\lambda \geq 0$ (respect.\ $\lambda \leq 0$). 
Lastly, by using Proposition \ref{prop:nobifz}(ii) again, there exists $\rho > 0$ small such that  $\Sigma_{0, \rho}^{\pm} := \Sigma^\pm_{0} \cap \overline{B_\rho ((0,0))}$ is closed and connected, and  if $(\lambda, u) \in \Sigma^\pm_{0, \rho} \setminus \{ (0,0) \}$ then $\lambda \gtrless 0$. So, $\mathcal{C}^\pm_0 := \Sigma^\pm_{0, \rho}$ have the desired properties.  \qed \newline

	%	% \begin{figure}[H] 

\subsection{Proof of assertion (ii) in Theorems \ref{resu01} and \ref{resu02}}
We consider the positivity of nontrivial nonnegative solutions of 
$(P_{\mathcal{B}})$ with $f(s)=s^q$, $q\in (0,1)$. 
Let $\mathfrak{S}$ be the nontrivial nonnegative solutions set of $(P_{\mathcal{B}})$, i.e.,
\begin{align*}
\mathfrak{S} :=  &  \{ (\lambda, u) \in\mathbb{R} \times V:0 \not \equiv u \geq0 \mbox{ solves } (P_{\mathcal{B}}) \}.
\end{align*}
Let $\mathfrak{C}$ be a nonempty connected subset of $\mathfrak{S}$, and let
\[
\mathfrak{C}^{\circ}:= \{ (\lambda, u) \in\mathfrak{C} : u \gg 0 \}.
\]

The following lemma (\cite[Theorem 1.7]{KRQU}) provides us with a nonexistence result for nontrivial nonnegative solutions of \eqref{prb:aqD}, 
which plays an important role in our argument when $\mathcal{B}u=\frac{\partial u}{\partial \mathbf{n}}$ and $\int_{\Omega}a \geq0$. 

\begin{lem}
\label{lem:nonon} Let $\mathcal{B}u=\frac{\partial u}{\partial \mathbf{n}}$. 
Assume $(H_\psi)$ with $\psi = a$. 
If $\int_{\Omega}a \geq0$ then there
exists $q_{a}^{*} \in(0,1)$ such that \eqref{prb:aqD} 
has no nontrivial nonnegative solutions for any $q \in(q_{a}^{*}, 1)$. 
\end{lem}

We give now sufficient conditions for the positivity of the nontrivial
nonnegative solutions on $\mathfrak{C}$ as follows. We recall 
that the sets $\mathcal{A}_{\mathcal{B}}^{\pm a}$ are given by \eqref{defA}.

\begin{prop}
\label{prop:posiD}
Let $f(s)=s^q$, $q\in (0,1)$. 
Suppose \eqref{Hg:posi}, $(H_\psi)$ with %\mpfns{modified}
$\psi=\pm a$, and the condition 
$b\geq 0$ in $\Omega$. Assume additionally $\int_\Omega a < 0$ if $\mathcal{B}u=\frac{\partial u}{\partial \mathbf{n}}$. If 
\begin{itemize}
\item $q \in \mathcal{A}_{\mathcal{B}}^{a}\cap\mathcal{A}_{\mathcal{B}}^{-a}$ (Dirichlet), 
\item $q\in \mathcal{A}_{\mathcal{B}}^{a} \cap \left( q_{-a}^{\ast},1 \right)$ (Neumann), 
\end{itemize}
where $q_{-a}^{\ast}$ is as in Lemma \ref{lem:nonon}, then 
$\mathfrak{C}^{\circ}$ is open and
closed in $\mathfrak{C}$. Consequently, $\mathfrak{C}^{\circ}= \mathfrak{C}$
if $\mathfrak{C}^{\circ}\neq\emptyset$.
\end{prop}

\noindent\textit{Proof.}
It is straightforward that $\mathfrak{C}^\circ$ is open in $\mathfrak{C}$, 
since $u \gg 0$ for $(\lambda, u) \in \mathfrak{C}^\circ$.
Next, we verify that $\mathfrak{C}^\circ$ is closed in $\mathfrak{C}$. Assume that $(\lambda_n, u_n) \in \mathfrak{C}^\circ$, and $(\lambda_n, u_n) \to (\lambda_0, u_0) \in \mathfrak{C}$ in $\mathbb{R} \times V$. 
We shall show that $(\lambda_0, u_0) \in \mathfrak{C}^\circ$. We discuss the following three cases, in accordance  with the sign of $\lambda_0$: 

\underline{(i) Case $\lambda_0 > 0$}. %\mpfns{added}
We use the condition $q\in \mathcal{A}_{\mathcal{B}}^{a}$ to deduce the desired assertion. In this case, $\lambda_n > 0$ for sufficiently large $n$. By the change of variables $v_n = \lambda_n^{-\frac{1}{1-q}}u_n$, we find that
\begin{align*}
\begin{cases}
-\Delta v_n = a(x)v_n^q + \lambda_n^{-\frac{1}{1-q}} b(x)g(\lambda_n^{\frac{1}{1-q}}v_n) & \mbox{ in } \Omega, \\
\mathcal{B}v_n = 0 & \mbox{ on } \partial \Omega.
\end{cases}
\end{align*}
Since $v_n \to v_0 = \lambda_0^{-\frac{1}{1-q}} u_0$ in $V$, we find that $v_0$ is a nonnegative weak solution of the problem
\begin{align*}
\begin{cases}
-\Delta v = a(x) v^q + \lambda_0^{-\frac{1}{1-q}} b(x) g(\lambda_0^{\frac{1}{1-q}}
v) & \mbox{ in } \Omega, \\
\mathcal{B}v = 0 & \mbox{ on } \partial \Omega.
\end{cases}
\end{align*}
In addition, $v_0 \not\equiv 0$ in $\Omega^a_+$. Indeed, since \eqref{Hg:posi} holds and $b\geq 0$, we see that
$v_n$ is a supersolution of \eqref{prb:aqD} which is positive in $\Omega$.
So, condition $(H_\psi)$ with $\psi = a$ allows us to apply \cite[Lemma 2.2]{KRQU}, and deduce that there exist a ball  $B \Subset \Omega_+^a$ and a continuous function $\psi$ on $\overline{B}$ such that $v_n \geq \psi > 0$ in $B$. Passing to the limit, we have
that $v_0\geq \psi$ in $B$, as desired.

By \eqref{Hg:posi} and the condition $b\geq 0$, $v_0$ is a supersolution of \eqref{prb:aqD}, and $v_0 > 0$ in $B$. On the other hand, we can construct a nonnegative subsolution $\psi_0$ of \eqref{prb:aqD} such that $\psi_0 \not\equiv 0$ in $B$, $\psi_0 \equiv 0$ in $\Omega \setminus B$, and $\psi_0 \leq v_0$. The sub and supersolutions method provides us with a solution $v_1$ of \eqref{prb:aqD} such that $\psi_0 \leq v_1 \leq v_0$, so that $v_1 \gg 0$, since $q \in \mathcal{A}_{\mathcal{B}}^{a}$. Consequently, we conclude that $u_0 = \lambda_0^{1/(1-q)} v_0 \gg 0$, as desired. 

\underline{(ii) Case $\lambda_0 < 0$}. 

(ii-1) Case $\mathcal{B}u=u$: %\mpfns{added}
We use the condition $q\in \mathcal{A}_{\mathcal{B}}^{-a}$ to deduce the desired assertion. In this case, $\lambda_n < 0$ for sufficiently large $n$. Setting $\mu = - \lambda$, $(P_{\mathcal{B}})$ turns into
\begin{align*}
\begin{cases}
-\Delta u = -\mu a(x) u^q + b(x)g(u) & \mbox{ in } \Omega, \\
u = 0 & \mbox{ on } \partial \Omega.
\end{cases}
\end{align*}
Set $\mu_n = - \lambda_n > 0$,  so that $v_n := \mu_n^{-\frac{1}{1-q}}u_n$ satisfies
\begin{align*}
\begin{cases}
-\Delta v_n = -a(x) v_n^q + \mu_n^{-\frac{1}{1-q}} 
b(x)g(\mu_n^{\frac{1}{1-q}}v_n) & \mbox{ in } \Omega, \\
v_n = 0 & \mbox{ on } \partial \Omega.
\end{cases}
\end{align*}
Since $v_n \to v_0 = \mu_0^{-\frac{1}{1-q}} u_0$ in $C^1_{0}(\overline{\Omega})$, we infer that
$v_0$ is a nonnegative weak solution of the problem
\begin{align*}
\begin{cases}
-\Delta v = -a(x) v^q + \mu_0^{-\frac{1}{1-q}} b(x) g(\mu_0^{\frac{1}{1-q}}v)
& \mbox{ in } \Omega, \\
v = 0 & \mbox{ on } \partial \Omega.
\end{cases}
\end{align*}
By using the condition $q \in \mathcal{A}_{\mathcal{B}}^{-a}$, the rest of the argument is carried out similarly as in the previous case.

(ii-2) Case $\mathcal{B}u=\frac{\partial u}{\partial \mathbf{n}}$: 
%\mpfns{added}
Under $q \in \left( q_{-a}^{\ast},1 \right)$, we shall see that case $\lambda_0 < 0$ does not occur, using Lemma \ref{lem:nonon}. 
We have $\lambda_n < 0$ for $n$ sufficiently large. By setting $\mu = - \lambda$, $(P_{\mathcal{B}})$ becomes 
\begin{align*}
\begin{cases}
-\Delta u = -\mu a(x) u^q + b(x) g(u) & \mbox{ in } \Omega, \\
\frac{\partial u}{\partial \mathbf{n}} = 0 & \mbox{ on } \partial \Omega.
\end{cases}
\end{align*}
Setting $\mu_n = - \lambda_n > 0$, and  $v_n = \mu_n^{-\frac{1}{1-q}}u_n$, we find that
\begin{align*}
\begin{cases}
-\Delta v_n = -a(x) v_n^q + \mu_n^{-\frac{1}{1-q}} b(x) g(\mu_n^{\frac{1}{1-q}}v_n) & \mbox{ in } \Omega, \\
\frac{\partial v_n}{\partial \mathbf{n}}  = 0 & \mbox{ on } \partial \Omega.
\end{cases}
\end{align*}
Since $v_n \to v_0 = \mu_0^{-\frac{1}{1-q}} u_0$ in $C^{1}(\overline{\Omega})$, we infer that $v_0$ is a nonnegative weak solution of the problem 
\begin{align*}
\begin{cases}
-\Delta v = -a(x) v^q + \mu_0^{-\frac{1}{1-q}} b(x) g(\mu_0^{\frac{1}{1-q}}v)
& \mbox{ in } \Omega, \\
\frac{\partial v}{\partial \mathbf{n}} = 0 & \mbox{ on } \partial \Omega.
\end{cases}
\end{align*}
In addition, $v_0 \not\equiv 0$ in $\Omega_{+}^{-a}$. Indeed, since 
\eqref{Hg:posi} holds, and $b\geq 0$, we see that  
$v_n$ is a positive supersolution of 
\begin{align} \label{prb:-aqN}
\begin{cases}
-\Delta v = -a(x) v^q & \mbox{ in } \Omega, \\
\frac{\partial v}{\partial \mathbf{n}} = 0 & \mbox{ on } \partial \Omega,
\end{cases}
\end{align} 
so that, by \cite[Lemma 2.2]{KRQU}, there exists a ball $B\Subset \Omega_{+}^{-a}$ and a continuous function $\psi$ on $\overline{B}$ such that $v_n \geq \psi > 0$ in $B$.
Passing to the limit, we obtain $v_0\geq \psi$ in $B$, as desired.

Now, we see that $v_0$ is also a nonnegative supersolution of \eqref{prb:-aqN} such that $v_0 > 0$ in $B$. Since we can construct a nonnegative subsolution $\psi_0$ of \eqref{prb:-aqN} such that $\psi_0 \not\equiv 0$ in $B$, $\psi_0 \equiv 0$ in $\Omega \setminus B$, and $\psi_0 \leq v_0$, the sub and supersolutions method provides us with a solution $v_1$ of \eqref{prb:-aqN} 
such that $\psi_0 \leq v_1 \leq v_0$ on $\overline{\Omega}$. So, $v_1$ is nontrivial and nonnegative. However, this contradicts Lemma \ref{lem:nonon}, since $\int_\Omega (-a) > 0$ and $q_{-a}^* < q < 1$. 

\underline{(iii) Case $\lambda_0 = 0$}. In this case, $u_0$ solves the problem
\[
\begin{cases}
-\Delta u_0 = b(x)g(u_0) & \mbox{ in } \Omega, \\
\mathcal{B}u_0 = 0 & \mbox{ on } \partial \Omega.
\end{cases}
\]
Since $u_0$ is nontrivial and nonnegative, the strong maximum principle and Hopf's lemma yield that $u_0 \gg 0$, as desired.

Lastly, since $\mathfrak{C}$ is connected, we conclude that $\mathfrak{C}^\circ = \mathfrak{C}$ if $\mathfrak{C}^\circ \neq \emptyset$. \qed \newline

Introducing the following growth condition on $g$: 
\begin{align} \label{H_2b'} 
0< \displaystyle{\lim_{s\to\infty} \frac{g(s)}{s^{p}}} < \infty 
\ \mbox{ for some } p>1, \ \mbox{where} \ p< \frac{N+1}{N-1} \ 
\mbox{if} \ N>2, 
\end{align} 
we can deduce that $\mathfrak{C}^\circ\neq \emptyset$, as shown (i) and (ii) below.

%\mpfnsb{joined two remarks, and modified}
\begin{rem} \label{rem:Cnon} 
\strut
%In (i) and (ii) below we give some sufficient conditions to have $\mathfrak{C}
%^{\circ}\neq\emptyset$ for $(P_{\mathcal{B}})$. 
\begin{enumerate}
\item[(i)] If $(0, u_{0}) \in \mathfrak{C}$ with $u_{0} \not \equiv 0$, 
then $u_{0}\gg 0$, i.e., $(0, u_{0}) \in\mathfrak{C}^{\circ}$. Indeed,
this is a direct application of the strong maximum principle and Hopf's lemma.

\item[(ii)] Assume \eqref{H_1g} and $(H_\psi)$ with 
$\psi = a$. Assume also \eqref{H_2'} if $\mathcal{B}u=u$;
%\marginpar{modified}\mpfns{agree}
and \eqref{H_2b'} and $\int_\Omega a < 0$ if $\mathcal{B}u=\frac{\partial u}{\partial \mathbf{n}}$. 
Let 
$b\equiv 1$ and $q \in\mathcal{A}_{\mathcal{B}}^{a}$. 
If there exist $(\lambda_{n}, u_{n}) \in\mathfrak{C}$
with $\lambda_{n} \to0^{+}$, then $u_{n} \gg 0$ for
sufficiently large $n$, i.e., $(\lambda_{n}, u_{n}) \in\mathfrak{C}^{\circ}$
for such $n$. This is an immediate consequence of \cite[Theorem 4.1, 
Theorem 4.5]{KRQU}.

\item[(iii)] When $\mathcal{B}u=\frac{\partial u}{\partial \mathbf{n}}$, 
Proposition \ref{prop:posiD} is valid for $\int_\Omega a > 0$, where 
we now assume $q \in \mathcal{A}_{\mathcal{B}}^{-a} \cap (q_{a}^{*}, 1)$ instead of $q \in \mathcal{A}_{\mathcal{B}}^{a} \cap (q_{-a}^{*}, 1)$. Indeed, when $\int_\Omega a > 0$, the case $\lambda_0 > 0$ does not occur, based on Lemma \ref{lem:nonon}, whereas the case $\lambda_0 < 0$ is verified as in the proof of Proposition \ref{prop:posiD}, using $\lambda a(x) = (-\lambda)(-a(x))$ and relying on \cite[Lemma 2.2]{KRQU}. 

\item[(iv)] Proposition \ref{prop:posiD} and %\mpfnsb{modified}
item (ii) hold more generally in the framework $\mathbb{R} \times C(\overline{\Omega})$, 
which can been seen by using elliptic regularity. \newline  
\end{enumerate}
\end{rem}

%We give the following further remarks for Proposition \ref{prop:posiD}.

%\begin{rem} \label{rem:inta>0} 
%\strut
%\begin{enumerate}
%\item When $\mathcal{B}u=\frac{\partial u}{\partial \mathbf{n}}$, 
%Proposition \ref{prop:posiD} is valid for $\int_\Omega a > 0$, where 
%we now assume $q \in \mathcal{A}_{\mathcal{B}}^{-a} \cap (q_{a}^{*}, 1)$ instead of $q \in \mathcal{A}_{\mathcal{B}}^{a} \cap (q_{-a}^{*}, 1)$. Indeed, when $\int_\Omega a > 0$, the case $\lambda_0 > 0$ does not occur, based on Lemma \ref{lem:nonon}, whereas the case $\lambda_0 < 0$ is verified as in the proof of Proposition \ref{prop:posiD}, using $\lambda a(x) = (-\lambda)(-a(x))$ and relying on \cite[Lemma 2.2]{KRQU}. 

%\item Proposition \ref{prop:posiD} and Remark \ref{rem:inta>0}(ii) 
%hold more generally in the framework $\mathbb{R} \times C(\overline{\Omega})$, 
%which can been seen by using elliptic regularity. \newline 
%\end{enumerate}
%\end{rem}

\noindent\textit{Proof of assertion (ii) in Theorem \ref{resu01}.}
We note from Remark \ref{rem:cc}(iii) that \eqref{Hg:posi} holds in case (b) of Theorem \ref{resu01}. 
Based on the result stated in Remark \ref{rem:Cnon}(i), this assertion is 
verified by a direct application of Proposition \ref{prop:posiD}. \qed\newline

\noindent\textit{Proof of assertion (ii) in Theorem \ref{resu02}.} 
Based on the result stated in Remark \ref{rem:Cnon}(ii), this assertion
is straightforward from Proposition \ref{prop:posiD} 
and Remark \ref{rem:Cnon}(iv). Indeed, we don't need to assume $(H_\psi)$ with $\psi = -a$ for applying Proposition \ref{prop:posiD} to the loop $\mathcal{C}_{\ast}$ given in Theorem \ref{resu02}(ii), since it lies in $\lambda \geq 0$ (see Figure \ref{fig:double}(ii)). 
Note that the condition $(H_\psi)$ with $\psi = -a$ is used only for case (ii) in the proof of Proposition \ref{prop:posiD}. \qed \newline

%\medskip
%\appendix 

%--------------------------------------------------------------------- 

\end{document}